\documentclass[12pt]{article}
\usepackage{fancyhdr, array,calc,graphicx,url,tabularx}
\usepackage{geometry}
\usepackage{latexsym}
\usepackage{amssymb}
\usepackage{amsmath}
\usepackage{makeidx}
\usepackage{enumerate}
\usepackage{verbatim}
\usepackage{tikz}
\usepackage[colorlinks=true,linkcolor=blue]{hyperref}
\usetikzlibrary{arrows.meta}
\allowdisplaybreaks

\usepackage{changepage}

\makeindex

{
   \end{minipage}
   \vspace*{\stretch{3}}
   \clearpage
}

\def\undertilde#1{{\baselineskip=0pt\vtop
  {\hbox{$#1$}\hbox{$\scriptscriptstyle\sim$}}}{}}

\newcommand{\utilde}{\undertilde}

\renewcommand{\models}{\vDash}

\newtheorem{theorem}{Theorem}[section]

\newtheorem{definition}[theorem]{Definition}

\newtheorem{lemma}[theorem]{Lemma}
\newtheorem{corollary}[theorem]{Corollary}
\newtheorem{claim}[theorem]{Claim}

\newtheorem{convention}[theorem]{Convention}
\numberwithin{figure}{section}

\newenvironment{Claim}{\noindent{\it
Claim} {}}{}
\newenvironment{Case}{\noindent{\it
Case}}{}

\def\a{\alpha}

\def\and{\mathrel{\kern1pt\&\kern1pt}}

\def\<#1>{\langle\,#1\,\rangle}

 \input xy
 \xyoption{all}

\title{Martin's Maximum${}^{++}$ implies Woodin's Axiom $(*)$}

\author{David Asper\'o\footnote{School of Mathematics, University of East Anglia, Norwich Research Park, Norwich, Norfolk, NR4 7TJ, UK. Funded by EPSRC Grant EP/N032160/1.} \ and Ralf Schindler\footnote{Institut f\"ur Mathematische Logik und Grundlagenforschung, Universit\"at M\"unster,
Einsteinstr. 62, 48149 M\"unster, FRG.
Funded by the Deutsche Forschungsgemeinschaft (DFG, German Research Foundation) under
Germany’s Excellence Strategy EXC 2044 390685587, Mathematics M\"unster: Dynamics - Geometry
- Structure.
}}   
         
\date{\today}

\DeclareMathOperator{\Ord}{Ord}

\DeclareMathOperator{\ZFC}{ZFC}
\DeclareMathOperator{\MA}{MA}
\DeclareMathOperator{\FA}{FA}

\DeclareMathOperator{\CH}{CH}

\DeclareMathOperator{\NS}{NS}

\begin{document}
   
\maketitle

\begin{abstract}
We show that Martin's Maximum${}^{++}$ implies Woodin's ${\mathbb P}_{\rm max}$ axiom $(*)$.
This answers a question from the 1990's and amalgamates two prominent axioms of set theory which were both known to imply that there are $\aleph_2$ many real numbers.
\end{abstract}

\section{Introduction.}

Cantor's Continuum Problem, which later became Hilbert's first Problem (see \cite{hilbert}), asks how many real numbers there are. After
having proved his celebrated theorem according to which ${\mathbb R}$ is uncountable, i.e., $2^{\aleph_0} > \aleph_0$, see \cite{cantor}, Cantor conjectured that every uncountable set of reals has the same size as ${\mathbb R}$, i.e., $2^{\aleph_0}=\aleph_1$. This statement is known as Cantor's Continuum Hypothesis ({\sf CH}). G\"{o}del \cite{goedel} proved in the 1930's that {\sf CH} is consistent with the standard axiom system for set theory, {\sf ZFC}, by showing that {\sf CH} holds in his constructible universe $L$, the minimal transitive model of {\sf ZFC} containing all the ordinals. The axiom $V=L$, saying that the universe $V$ of all sets is simply identical with $L$, has often been rejected, however, as an undesirable minimalistic assumption about $V$. For instance, $L$ cannot have measurable cardinals by a result of Scott \cite{scott}. G\"{o}del himself believed that {\sf CH} would be shown not to follow from {\sf ZFC}, and at least for part of his life he held the view that {\sf CH} is indeed false and that actually 
\begin{eqnarray}\label{-ch}
2^{\aleph_0}=\aleph_2
\end{eqnarray}
(see \cite{godel_nonCH} and
\cite[pp.\ 173ff.]{goedel-feferman}). 
In 1947, G\"odel \cite{what_is_ch} wrote:

\begin{quote}
[...] one may on good reason suspect that the role of the continuum problem in set theory will be this, that it will finally lead to the discovery of new axioms which will make it possible to disprove Cantor's conjecture.
\end{quote}

As we shall now try to explain, in the light of our unifying result, Theorem \ref{main_thm}, one could make the case that
with the two axioms ${\sf MM}^{++}$ and $(*)$, natural and strong such axioms
have already been found.

Luzin \cite{luzin} proposed a related hypothesis which also refutes {\sf CH}, namely 
\begin{eqnarray}\label{luzin}
2^{\aleph_0}=2^{\aleph_1}.
\end{eqnarray}

That {\sf CH} does not follow from {\sf ZFC} was confirmed by Cohen in 1963 through the discovery of the method of {\em forcing}: Every model of {\sf ZFC} can be generically extended to a model of {\sf ZFC} in which {\sf CH} fails, see \cite{cohen}. In fact, using forcing one can show that it is relatively consistent with {\sf ZFC} that the cardinality of the continuum is $\aleph_1$, $\aleph_2$, $\aleph_{155}$, $\aleph_{\omega^2 + 17}$, or  $\aleph_\alpha$ for many other values of $\alpha$, see \cite{solovay}.  

\subsection{New axioms.} 

Ever since Cohen's work, set-theorists have been searching for natural new axioms which extend {\sf ZFC} and which settle the Continuum Problem (see e.g.\ \cite{hugh-ch1}, \cite{hugh-ch2}, \cite{peter}, and the discussion in
\cite{new?}). One family of such axioms is the hierarchy of large cardinal axioms. It was realized early on, though, that these axioms cannot settle the Continuum Problem: one can always force {\sf CH} to hold or be false by small forcing notions, and all large cardinals which exist in $V$ will retain their large cardinal properties in the respective extensions, see \cite{l-s}. 

One axiom which does settle the Continuum Problem is {\sf CH} itself; after all, {\sf CH} looks natural in that it gives the least possible value to $2^{\aleph_0}$ consistent with Cantor's theorem, $2^{\aleph_0}>\aleph_0$. {\sf CH} allows the ``diagonal'' construction of objects of size $\aleph_1$ with specific combinatorial properties, e.g.\ Luzin or Sierpi\'nski sets.
In 1985, Woodin proved his $\Sigma^2_1$ absoluteness result conditioned on {\sf CH}.
Namely, if {\sf CH} holds true, there is a proper class of measurable Woodin cardinals, and $\sigma$ is a statement of the form ``There is a set of reals $X$ such that $\varphi(X, r)$,'' where $r$ is a real and $\varphi(X)$ is a formula of set theory all of whose quantifiers are restricted to 
reals, such that $\sigma$ can be forced over $V$, then $\sigma$ actually holds true in $V$, see e.g.\ \cite[Theorem 4.1]{ds2}. 
Over the last decade, Woodin has developed a sophisticated scenario for set theory according to which {\sf CH} is true, see e.g.\ \cite{ICM} and \cite{midrasha}.

Nevertheless, and despite the appeal of $\Sigma^2_1$ absoluteness, {\sf CH} is often regarded as 
a minimalistic assumption on a par with its parent, $V=L$.
To give an illustrative example,
under {\sf CH} one can easily find sets $X$ and $Y$ of reals without endpoints which are both $\aleph_1$-dense---in the sense that every interval of points 
contains exactly $\aleph_1$ many points---but such that $X$ and $Y$ are not order-isomorphic. On the other hand,
by a theorem of Baumgartner \cite{baumgartner}, given any such $X$ and $Y$, there is a
nicely behaved forcing notion which adds an order-isomorphism between $X$ and $Y$. Thus, adopting {\sf CH} 
precludes the existence of sufficiently generic filters for such forcing notions---which may consistently 
exist.

A dual approach to {\sf CH} is to formulate axioms stipulating the existence of objects which may possibly exist, i.e., 
to look for ``maximality principles'' expressing some form of saturation of the universe of all sets with respect to its generic extensions. Such principles are known as \emph{forcing axioms}. 
Shortly after the discovery of forcing, it was realized that it is possible to iterate the process of forming generic extensions $V \subset V[g_0] \subset V[g_1] \subset \ldots \subset V[g_\alpha]$ of $V$ in any length $\alpha$ in such a way that the final model is itself a generic extension of $V$.
By ``closing off'' one may then get to final models which are in fact saturated with respect to the existence of certain (partial) generics in the way prescribed by forcing axioms.

\subsection{Forcing axioms.}

Forcing axioms are generalizations of the Baire Category Theorem. Formally, they assert the existence of sufficiently generic filters for all members of some reasonably large class of forcing notions. In a general form, given an infinite cardinal $\kappa$ and a class $\mathcal K$ of forcing notions, the forcing axiom $\FA_\kappa(\mathcal K)$ is the statement that for every $\mathbb P\in\mathcal K$ and for every collection $\mathcal D$ of dense subsets of $\mathbb P$ such that $|\mathcal D|=\kappa$ there is a filter $g$ of $\mathbb P$  which is $\mathcal D$-generic (i.e., is such that $g \cap D \neq \emptyset$ for each $D\in\mathcal D$).  $\FA_\kappa(\mathcal K)$ is to be seen as a maximality principle with respect to forceability via forcing notions from $\mathcal K$: If $\FA_\kappa(\mathcal K)$ holds, then all $\Sigma_1$ statements with parameters in $H_{\kappa^+}$ that can be forced to hold by some forcing notion in $\mathcal K$ already hold in the universe.\footnote{As a matter of fact, forcing axioms of the form $\FA_\kappa(\mathcal K)$ can 
be fully characterized in terms of a suitable form of $\Sigma_1$-absoluteness with respect to generic extensions via members from $\mathcal K$ (see e.g.\ \cite{hugh} or \cite[Theorem 1.3]{cs2}).} The answers to questions about $H_{\kappa^+}$ provided by forcing axioms $\FA_\kappa(\mathcal K)$ are often regarded as being natural in that $\FA_{\kappa^+}(\mathcal K)$ offers a uniformly `saturated' picture of $H_{\kappa^+}$, ruling out the type of pathological objects that one can construct when $H_{\kappa^+}$ has an artificially constrained structure. 

In what follows we will consider only forcing axioms $\FA_\kappa(\mathcal K)$ for $\kappa=\omega_1$.\footnote{This is the first $\kappa$ for which $\FA_\kappa(\mathcal K)$ does not follow outright from $\ZFC$. Also, $\kappa=\omega_1$ is the only level for which we currently have a reasonably complete picture of the available forcing axioms.} The first such forcing axiom shown to be consistent was Martin's Axiom at $\omega_1$, $\MA_{\omega_1}$, see \cite{solovay-tennenbaum} and \cite{martin-solovay}.  $\MA_{\omega_1}$ is $\FA_{\omega_1}(\mathcal K)$, where $\mathcal K$ is the class of partial orders $\mathbb P$ with the countable chain condition (i.e., such that there is no uncountable family of pairwise incompatible conditions in $\mathbb P$). 
Over the following years, a number of generalizations of $\MA_{\omega_1}$ were isolated. {\sf PFA}, the Proper Forcing Axiom, is 
$\FA_{\omega_1}(\mathcal K)$, where $\mathcal K$ is the class of partial orders $\mathbb P$ which are proper.

The following is a list of examples of natural statement which are implied by forcing axioms.\label{list!} 

\begin{itemize}

\item $\MA_{\omega_1}$ implies that there are no Suslin lines (\cite{solovay-tennenbaum}).
\item $\MA_{\omega_1}$ implies that every union of $\aleph_1$-many Lebesgue null subsets of reals is Lebesgue null (\cite{martin-solovay}). 
\item $\MA_{\omega_1}$ implies the existence of a non-free Whitehead group (\cite{Shelah-Whitehead}). 
\item {\sf PFA} implies Baumgartner's Axiom that all $\aleph_1$-dense sets of reals are order-isomorphic (essentially \cite{baumgartner}).\footnote{But it does not follow from $\MA_{\omega_1}$ (s.\ \cite{ARS}).}
\item {\sf PFA} implies Kaplansky's conjecture (\cite{todo}).
\item {\sf PFA} implies that there is a $5$-element basis for the class of uncountable linear orders (\cite{Moore5}). 
\item {\sf PFA} implies that every automorphism of the Calkin algebra of a separable Hilbert space is inner (\cite{Farah}). 
\end{itemize}



This line of research culminated in the proof by Foreman-Magidor-Shelah of the consistency of \emph{Martin's Maximum}, {\sf MM}, see \cite{FMS}. 

Martin's Maximum is $\FA_{\omega_1}(\mathcal K)$, where $\mathcal K$ is the class of partial orders $\mathbb P$ such that forcing with $\mathbb P$ preserves the stationarity of all stationary subsets of $\omega_1$ in $V$. {\sf MM} is provably maximal in the sense that the forcing axiom  $\FA_{\omega_1}(\{\mathbb P\})$ fails for any forcing notion $\mathbb P$ destroying some stationary subset $S$ of $\omega_1$.
At the same time, {\sf MM} can be forced by means of a forcing iteration ${\mathbb P} \subset V_\kappa$, assuming that $\kappa$ is a supercompact cardinal. 
The natural such forcing ${\mathbb P}$ actually produces a model of a strengthening of {\sf MM}, called ${\sf MM}^{++}$. This is the statement that if $\mathbb P$ is a forcing notion preserving stationary subsets of $\omega_1$,  $\mathcal D$ is a collection of size $\aleph_1$ consisting of dense subsets of $\mathbb P$, and $\{\tau_\alpha\,:\,\alpha<\omega_1\}$ is a collection of $\mathbb P$-names for stationary subsets of $\omega_1$, then there is a filter $g \subset \mathbb P$ which is $\mathcal D$-generic and which, furthermore, interprets every $\tau_\alpha$, $\alpha<\omega_1$, as a truly stationary set in $V$ (i.e., $\{\nu<\omega_1\,:\, \exists p \in g, \, p\Vdash_{\mathbb P} {\check \nu} \in\tau_\alpha \}$ is stationary for every $\alpha<\omega_1$). 

Already $\MA_{\omega_1}$ contradicts {\sf CH}, and it even proves Luzin's hypothesis (\ref{luzin}), i.e., $2^{\aleph_0}=2^{\aleph_1}$. More interestingly, {\sf MM} (in contrast to $\MA_{\omega_1}$) decides the cardinality of the continuum,
and in fact it 
confirms G\"odel's conjecture
(\ref{-ch}), $2^{\aleph_0}=\aleph_2$. This is shown by producing an affirmative answer to Friedman's Problem under {\sf MM}, see \cite[Theorems 9 and 10]{FMS}.\footnote{It was later verified by Moore that already the much weaker forcing axiom {\sf BPFA} implies (\ref{-ch}), see \cite{justin}.} ${\sf MM}^{++}$ is---by its very definition and the fact that no strictly stronger forcing axiom can be consistent---a prototype maximality principle for $V$. Remarkably, the empirical evidence seems to suggest that ${\sf MM}^{++}$ provides a complete theory of the initial segment $H_{\omega_2}$ of the universe of sets, at least with respect to natural questions. 
Here, $H_{\omega_2}$ is the collection of all sets which are hereditarily of size $< \aleph_2$.

\subsection{The ${\mathbb P}_{\rm max}$ axiom $(*)$.}

There is another maximality principle, though, which Magidor called a ``competitor'' of {\sf MM}, see \cite[p.\ 18]{menachem}, and which is denoted by $(*)$. Its formulation involves the notion of ${\mathbb P}_{\rm max}$, a forcing which
was isolated by W.H.\ Woodin, see \cite[Definition 4.33]{hugh} and Definition \ref{defn_P_max} below. In much the same way as {\sf MM}, $(*)$ is inspired by and formulated in the language of forcing, and they both have ``the same intuitive motivation: Namely, the
universe of sets is rich'' (\cite[p.\ 18]{menachem}). 
$(*)$, introduced by Woodin in \cite[Definition 5.1]{hugh}, is the conjunction of the following two statements.
\begin{enumerate}
\item[(i)]
{\sf AD}, the Axiom of Determinacy,\footnote{See e.g.\ \cite[Chapter 12]{book}.} holds in $L({\mathbb R})$, and 
\item[(ii)] there is some $g$
which is ${\mathbb P}_{\rm max}$-generic over $L({\mathbb R})$ such that
${\cal P}(\omega_1) \subset L({\mathbb R})[g]$.
\end{enumerate}

Item (i), that ${\sf AD}$ holds in ${L(\mathbb R)}$, follows from the existence of large cardinals, e.g.\ from the existence of infinitely many Woodin cardinals with a measurable cardinal above them all, see \cite[p.\ 91]{proj-determinacy}. Item (ii) is the part of $(*)$ which goes beyond assuming the existence of large cardinals.  $\mathbb P_{\max}$ 
arose out of earlier work by Steel-Van Wesep \cite{StVW} and by Woodin \cite{hugh83}
on the size of
\boldmath$\delta$\unboldmath${}^1_2$ and the question if ${\sf NS}_{\omega_1}$, the nonstationary ideal on $\omega_1$, can be saturated.
Here, ${\sf NS}_{\omega_1}$ is called saturated iff there is no collection ${\cal A}$ of stationary subsets of $\omega_1$ of size $\aleph_2$ such that $S \cap T$ is nonstationary for all $S$, $T \in {\cal A}$, $S \not= T$.

${\mathbb P}_{\rm max}$ consists of countable transitive structures, membership in ${\mathbb P}_{\rm max}$ is uniformly $\Pi^1_2$ in the codes, and the order $<_{{\mathbb P}_{\rm max}}$ is arithmetical. ${\mathbb P}_{\rm max}$ is $\omega$-closed and homogeneous, see \cite[Lemma 4.43]{hugh}. 
The fact that a forcing ${\mathbb P}$ is homogeneous means that the validity in the forcing extension of a given statement is decided in the ground model by the trivial condition in ${\mathbb P}$.
The homogeneity of ${\mathbb P}_{\rm max}$ then yields that under $(*)$, the theory of $L({\mathcal P}(\omega_1))$ becomes part of the theory of $L({\mathbb R})$ in the sense that if $\varphi$ is any sentence, then 
\begin{eqnarray}\label{part-of}
L({\mathcal P}(\omega_1)) \models \varphi \mbox{ if and only if } \Vdash_{L({\mathbb R})}^{{\mathbb P}_{\rm \max}} \varphi.
\end{eqnarray} 
If {\sf AD} holds in ${L(\mathbb R)}$, then there is no well-order of the reals in ${L(\mathbb R)}$ (see 
e.g.\ \cite[Lemma 12.2]{book}), but if $g$ is ${\mathbb P}_{\rm max}$-generic over ${L(\mathbb R)}$, then {\sf ZFC} is true in $L({\mathbb R})[g]$ (see \cite[Theorem 4.54]{hugh}), and moreover ${\sf NS}_{\omega_1}$ is saturated in $L({\mathbb R})[g]$ (see \cite[Theorem 4.50]{hugh}) 
and $L({\mathbb R})[g]$ provides an effective failure of {\sf CH} in that \boldmath$\delta$\unboldmath${}^1_2 = \omega_2$ is true in $L({\mathbb R})[g]$ (see \cite[Theorem 4.53]{hugh}).

Like the ``classical'' forcing axioms culminating with ${\sf MM}^{++}$, $(*)$ is also a maximality principle. While $(*)$ implies none of the stronger forcing axioms, see e.g.\ \cite[Theorem 1.3]{mlq0002d}, it does imply ${\sf MA}_{\omega_1}$. In particular, $(*)$ implies the first three implications of ${\sf MA}_{\omega_1}$ which are listed on p.\ \pageref{list!}. As it turns out,  $(*)$ also implies that every automorphism of the Calkin algebra of a separable Hilbert space is inner, see \cite{paul!}, \cite{Farah}; and, at least in conjunction with the existence of a Woodin cardinal,\footnote{This is by the proof of Theorem \ref{Wo1}.} it also implies Baumgartner's Axiom on $\aleph_1$-dense sets of reals, as this can be expressed by a $\Pi_2$ sentence over $H_{\omega_2}$, as well as the existence of a $5$-element basis for the uncountable linear orders (since, in the presence of Baumgartner's Axiom and $\MA_{\omega_1}$, the existence of such a basis follows from every Aronszajn line containing a Countryman line, see \cite{Moore5}, which again can be expressed by a  $\Pi_2$ sentence over $H_{\omega_2}$). 

$(*)$ implies (and is in fact 
equivalent to) what is dubbed ``$\Pi_2$ maximality.'' 
A sentence $\sigma$ (in the language of set theory, possibly augmented with some additional predicates) is said to be $\Pi_2$ if it is of the form $\forall x\exists y\varphi(x, y)$, with $\varphi(x, y)$ being a formula with only restricted quantifiers. 
There is a whole 
family of interesting statements which are $\Pi_2$ in the language for the structure
$$(H_{\omega_2}; \in, \NS_{\omega_1}){\rm , }$$
see e.g.\ the discussion in \cite{ds}.
The formulation of ``$\Pi_2$ maximality'' involves the concept of $\Omega$-logic, see 
\cite[Section 10.4]{hugh}; for a sentence $\sigma$ to be ``$\Omega$-consistent'' is stronger than it just being consistent in that $\sigma$ needs to be true in models which are closed under arbitrarily complicated universally Baire operations, see \cite[Definition 10.144]{hugh}. The $\Pi_2$ maximality theorem, see \cite[Theorem 10.150]{hugh}, then runs as follows.

\begin{theorem}\label{Wo1} {\em \textbf{(Woodin)}} Suppose there is a proper class of Woodin cardinals. Then the following statements are equivalent.

\begin{enumerate}
\item[(1)] $(*)$.
\item[(2)] 
Let $\sigma$ be a $\Pi_2$ sentence in the language for the structure
$$(H_{\omega_2}; \in, \NS_{\omega_1}, A \colon A \in {\mathcal P}({\mathbb R}) \cap L({\mathbb R})).$$ 
If $\sigma$ is $\Omega$-consistent, then $\sigma$ is true.
\end{enumerate}
\end{theorem}

One specific instance of $\sigma$ in (2) of Theorem \ref{Wo1} is called $\psi_{\rm AC}$, see \cite[Definition 5.12]{hugh}. It is in spirit a
local version of an affirmative solution to Friedman's Problem. Woodin showed, see \cite[Theorem 5.14, Lemmata 5.15 and 5.18]{hugh}, that $\psi_{\rm AC}$
follows from both {\sf MM} and $(*)$ and that $\psi_{\rm AC}$ implies
G\"odel's conjecture (\ref{-ch}), i.e.\ $2^{\aleph_0}=\aleph_2$. 

The homogeneity of ${\mathbb P}_{\rm max}$ gives that 
in the presence of large cardinals, 
$(*)$ yields a complete theory for $L({\cal P}(\omega_1))$ modulo set-forcing: by (\ref{part-of}), all set-generic extensions of $V$ in which $(*)$ holds true agree on the theory of $L({\cal P}(\omega_1))$.

Despite its nice properties, in order for $(*)$ to be a convincing candidate for a natural axiom, it would 
have to be compatible with all consistent large cardinal axioms. While $L(\mathbb R)[g]$
is trivially a model of
$(*)$, provided that $g$ is ${\mathbb P}_{\rm max}$-generic over $L(\mathbb R)$, 
Scott's result \cite{scott} carries over from $L$ to
$L(\mathbb R)[g]$ and shows that $L(\mathbb R)[g]$ cannot have measurable cardinals either.
\cite{hugh} and subsequent work left open the problem whether $(*)$ would be compatible with large cardinals beyond the level of Woodin cardinals. 

\subsection{Unifying forcing axioms and $(*)$.}

Prior to the current paper, the relation between classical forcing axioms like {\sf MM}, 
which could be forced by iterated forcing over models of {\sf ZFC} with large cardinals, and the axiom $(*)$, whose models were obtained by forcing over models satisfying 
the Axiom of Determinacy, 
remained a complete mystery. It had been known by a result of P.\ Larson \cite{plarson} that even ${\sf MM}^{+\omega}$, an axiom strictly between {\sf MM} and ${\sf MM}^{++}$, does not imply $(*)$.\footnote{${\sf MM}^{+\omega}$ is the strengthening of ${\sf MM}$ obtained by replacing, in the formulation of ${\sf MM}^{++}$, collections of $\aleph_1$ many names for stationary sets with collections of only countably many such names.} One can build models of ${\sf MM}^{+\omega}$ with a well-order of $H_{\omega_2}$ which is definable over $(H_{\omega_2}; \in)$ by a formula without parameters, and the existence of such a well-order is incompatible with $(*)$ by the homogeneity of ${\mathbb P}_{\rm max}$. It remained even unclear whether classical strong forcing axioms would be compatible at all with $(*)$, see 
\cite[p.\ 846]{hugh}. See also \cite[Question (18) a) on p.\ 924]{hugh}, \cite[Conjecture 6.8 on p.\ 19]{menachem},
and \cite[Problem 14.7]{moore}.

This paper resolves the tension between {\sf MM} and $(*)$. We prove:

\begin{theorem}\label{main_thm}
Assume {\em Martin's Maximum}${}^{++}$. Then Woodin's ${\mathbb P}_{\rm max}$-axiom $(*)$ holds true.
\end{theorem}

In particular, {\sf MM} and $(*)$ are compatible with one another, and
$(*)$ is compatible with all consistent large cardinal axioms: If $\kappa$ is a supercompact cardinal and ${\mathbb P} \subset V_\kappa$ is the partial order from
\cite{FMS}
to force
${\sf MM}^{++}$, then by Theorem \ref{main_thm} the axiom $(*)$ holds in $V^{\mathbb P}$, and all the large cardinals of $V$ above $\kappa$ are preserved by ${\mathbb P}$. 

Theorem \ref{main_thm} renders ${\sf MM}^{++}$ a particularly appealing axiom. Not only is ${\sf MM}^{++}$ a provably maximal forcing axiom providing the `right' answers to questions pertaining to $H_{\omega_2}$,\footnote{The `right' answers from a conception of the universe as being uniformly saturated with respect to forcing. 
} but it follows from Theorem \ref{main_thm} that ${\sf MM}^{++}$  implies the form of $\Pi_2$ maximality for arbitrary set-forcing given by (2) of Theorem \ref{Wo1} and, moreover, that ${\sf MM}^{++}$ completely decides the theory of $L(\mathcal P(\omega_1))$ via set-forcing. 

It also follows from Theorem \ref{main_thm} 
that $(*)$ can be characterized, in the presence of large cardinals, by a statement which on the face of its formulation is weaker than (2) of Theorem \ref{Wo1}. We will prove the following theorem at the end of the next section.

\begin{theorem}\label{char_of_star} Suppose there is a supercompact cardinal. Then the following statements are equivalent.

\begin{enumerate}
\item[(1)] $(*)$.
\item[(2)] 
Let $\sigma$ be a $\Pi_2$ sentence in the language for the structure
$$(H_{\omega_2}; \in, \NS_{\omega_1}, A \colon A \in {\mathcal P}({\mathbb R}) \cap L({\mathbb R})).$$ 
If there is a stationary set preserving forcing ${\mathbb P}$ such that $\sigma$
holds in $V^{\mathbb P}$, then $\sigma$ is true in $V$.
\end{enumerate}
\end{theorem}

This equivalence of $(*)$ is a variant of one which we are going to state below, see Theorem
\ref{char_star}, and which characterizes $(*)$ as a strong version of a bounded forcing axiom.

The authors thank Ilijas Farah, Andreas Lietz, and Matteo Viale for their comments on earlier versions of this paper. They would also like to thank the anonymous referees for their reports,
which have been 
of extraordinary quality and most helpful. Finally, they would like to thank Andreas Lietz for drawing the beautiful diagrams. 

The reader should have some acquaintance with forcing, determinacy, and universally Baire sets of reals. The relevant material is covered e.g.\ in \cite[Chap.\ 6 and sections 7.1, 8.1, and 12.1]{book}. Familiarity with stationary set preserving forcings and Martin's Maximum, see e.g.\ \cite{FMS} or \cite[Chap.\ 37]{jech}, and with ${\mathbb P}_{\rm max}$ forcing to the extent of say \cite[Chap.\ 4]{hugh} or \cite[Sections 1-6]{larson-handbook} would be desirable. Knowledge of forcings which are similar to the one which will be designed here and which were developed earlier e.g.\ in \cite{ronald_l_forcing}, \cite{cs}, \cite{ds}, or \cite{ds2} is by no means required or presupposed. 

\section{Preliminaries.}\label{section2}

Let us first state again 
{\em Martin's Maximum}${}^{++}$, abbreviated by {\sf MM}${}^{++}$ and isolated by Foreman-Magidor-Shelah \cite{FMS} (cf.\ also \cite[Definition 2.45 (2)]{hugh}). 

\begin{definition}
{\sf MM}${}^{++}$ is the statement that if ${\mathbb P}$ is a forcing which preserves stationary subsets of $\omega_1$, if $\{ D_i \colon
i<\omega_1 \}$ is a collection of dense subsets of ${\mathbb P}$, and if $\{ \tau_i
\colon i<\omega_1 \}$ is a collection of ${\mathbb P}$-names for stationary subsets
of $\omega_1$, then there is a filter $g \subset {\mathbb P}$ such that 
for every $i<\omega_1$,
\begin{enumerate}
\item[(i)]
$g \cap
D_i \not= \emptyset$ and 
\item[(ii)] $(\tau_i)^g = \{ \xi<\omega_1
\colon \exists p \in g \, p \Vdash_{\mathbb P} {\check \xi} \in \tau_i \}$ is stationary.
\end{enumerate}
\end{definition}

The forcing ${\mathbb P}_{\rm max}$ was designed by W.\ Hugh Woodin, see \cite[Chapter 4]{hugh}, specifically \cite[Definition 4.33]{hugh}. 

In order to define ${\mathbb P}_{\rm max}$, we need the notion of
``generic iterability'' of structues of the form $(M; \in, I,a)$, where $M$ is a transitive model of a sufficiently large fragment of ${\sf ZFC}$, $(M;I)$ is amenable,\footnote{I.e., $x \cap I \in M$ for all $x \in M$.} $a \subset \omega_1^M$,
and $(M;I) \models$ ``$I$ is a normal uniform ideal on $\omega_1$.'' Given an ordinal $\gamma\leq\omega_1$, $\langle \langle (M_i; \in, I_i,a_i)
\colon i \leq \gamma \rangle, \langle \pi_{i,j} \colon i \leq j \leq \gamma \rangle , \langle g_i \colon i < \gamma \rangle \rangle$ is a generic iteration of $(M; \in, I, a)$ if the following hold true.
\begin{itemize}
\item $(M_0;\in,I_0,a_0)=(M;\in,I,a)$,
\item for $i<\gamma$, $g_i$ is a $\mathcal P(\omega_1)^{M_\alpha}\setminus I_i$-generic filter over $M_i$, $M_{i+1}$ is the ultrapower of $M_i$ by $g_i$, and $\pi_{i,i+1} \colon (M_i; \in, I_i,a_i)\to (M_{i+1}; \in, I_{i+1},a_{i+1})$ is the corresponding generic elementary embedding, 
\item $\pi_{i, k}=\pi_{j,k}\circ \pi_{i,j}$ for all $i \leq j \leq k$, and
\item if $\beta$ is a nonzero limit ordinal $\leq \gamma$, then $(M_\beta , (\pi_{i, \beta} \colon i<\beta))$ is the direct limit of $(M_i , \pi_{i, j}
\colon i \leq j<\beta)$. 
\end{itemize}

\noindent $(M;\in,I,a)$ being generically iterable means that all models in any generic iteration of $(M; \in, I,a)$ are well-founded, irrespective of the filters $g_\alpha$ chosen at any stage $\alpha$, 
see \cite[Definition 4.1]{hugh}. 

Let us stress that the current paper will only consider such generic iterations rather than iterations of mice as being studied in inner model theory. 

\begin{definition}\label{defn_P_max}
The conditions in ${\mathbb P}_{\rm max}$ are countable transitive models of a sufficiently large fragment of ${\sf ZFC}$ plus ${\sf MA}_{\omega_1}$ of the form $(M;\in,I, a)$, where 
\begin{enumerate}
\item[(i)] $(M;I)$ is amenable
and $(M;I) \models$ ``$I$ is a normal uniform ideal on $\omega_1$,''
\item[(ii)] $a \in {\cal P}(\omega_1^M) \cap M$ and $M \models$ ``$\omega_1 =
\omega_1^{L[a,x]}$ for some real $x$,'' and
\item[(iii)] $(M;\in,I)$ is generically iterable.
\end{enumerate}
We construe ${\mathbb P}_{\rm max}$ as a partial order by declaring that $(N;\in,J,b)$ is stronger than $(M;\in,I,a)$, denoted by $(N;\in,J,b) < (M;\in,I,a)$, if and only if $(M;\in,I,a) \in N$ and inside $N$ there is a generic iteration of $(M;\in,I,a)$ of length $\omega_1^N+1$ with last model $(M^*;\in,I^*,a^*)$ such that $I^* = J \cap M^*$ and $a^* = b$.\footnote{${\mathbb P}_{\rm max}$ according to our Definition \ref{defn_P_max} is a slightly bigger poset than the one according to 
\cite[Definition 4.33]{hugh}. The difference is that we weakened the requirement $I \in M$ of \cite[Definition 4.33]{hugh} to ``$(M;I)$ is amenable.'' This natural move will make $(H_{\omega_2};\in,{\sf NS}_{\omega_1},A)$ for any $A \subset \omega_1$ a ${\mathbb P}_{\rm max}$ condition in a generic extension where $H_{\omega_2}$ is countable, cf.\ (\ref{its-a-Pmax-condition!}). It is easy to see that ${\mathbb P}_{\rm max}$ according to \cite[Definition 4.33]{hugh} is dense in ${\mathbb P}_{\rm max}$ according to our Definition \ref{defn_P_max}, so that both forcing notions are forcing-equivalent.}
\end{definition}

Most of \cite{hugh} studies the effect of forcing with ${\mathbb P}_{\rm max}$ or variants thereof over a model of the Axiom of Determinacy.
Let us state again 
Woodin's ${\mathbb P}_{\rm max}$ axiom $(*)$, see \cite[Definition 5.1]{hugh}. 

\begin{definition}
$(*)$ says that 
\begin{enumerate}
\item[(i)]
{\sf AD} holds in $L({\mathbb R})$ and 
\item[(ii)] there is some $g$
which is ${\mathbb P}_{\rm max}$-generic over $L({\mathbb R})$ such that
${\cal P}(\omega_1) \subset L({\mathbb R})[g]$.
\end{enumerate}
\end{definition}

Already the {\em Proper Forcing Axiom}, {\sf PFA}, which is much weaker than {\sf MM}${}^{++}$, implies ${\sf AD}^{L({\mathbb R})}$ and much more, see \cite{john}, \cite{stacking}, and \cite[Chapter 12]{trang-sargsyan}.

The current paper produces a proof of Theorem \ref{main_thm}.
%
Our key new idea is ($\Sigma$.8) on page \pageref{(C.9)} below.
We try to give an overview of the proof of Theorem \ref{main_thm} at the end of this section.


Theorem \ref{main_thm} is optimal in that P.\ Larson \cite{plarson} and \cite{plarson2} has shown that ${\sf MM}^{+\omega}$ is consistent with $\lnot (*)$ 
relative to a supercompact limit of supercompact
cardinals.
Our proof is also 
optimal
in that the forcing which we will use to verify Theorem \ref{main_thm} has size 
$2^{\aleph_2}$, while Woodin has shown that ${\sf MM}^{++}$ for forcings 
of size $2^{\aleph_0} = \aleph_2$ does not imply $(*)$, see \cite[Theorem 10.90]{hugh}, and it is consistent with ${\sf MM}^{++}$ that $2^{\aleph_2}=
\aleph_3$.

Throughout our entire paper, ``$\omega_1$'' 
will {\em always} denote $\omega_1^V$, the $\omega_1$ of $V$. We shall also make permanent use of the following.  

\begin{convention}\label{convention}
Let us fix throughout this paper some $A \subset \omega_1$ such that $\omega_1^{L[A]} = \omega_1$. Let us 
define $g_A$ as the set of all ${\mathbb P}_{\rm max}$ conditions $p=(N;\in,I,a)$
such that there is a generic iteration $$(N_i,\sigma_{ij} \colon i \leq j \leq \omega_1)$$
of $p = N_0$ of length $\omega_1+1$ such that if we write $N_{\omega_1} = (N_{\omega_1};\in,I^*, a^*)$,\footnote{Here and elsewhere we often confuse a model with its underlying universe.} then $I^* = ({\sf NS}_{\omega_1})^V \cap N_{\omega_1}$ and $a^* = A$.
\end{convention}

In the following statement, $X^\#$ denotes the sharp of $X$. While the formal definition of a sharp (see e.g. \cite[Section 10.2]{book}) won't play any role in what follows, the reader may think of ``${\cal P}(\omega_1)^\#$ exists'' as just some extra large cardinal structure which is assumed to be present in the universe.

We are going to use now the concept of elementary substructures.
For any two models ${\cal M}$ and ${\cal N}$ with underlying universes $M$ and $N$, respectively, and with the same first order language associated to them,
${\cal M} \prec {\cal N}$ means that ${\cal M}$ is an elementary substructure of ${\cal N}$, i.e., $M \subset N$ and for all formulae $\varphi$ of that common language and
all $x_1$, $\ldots$, $x_k \in M$, 
\begin{eqnarray}\label{ee}
{\cal M} \models \varphi(x_1,\ldots,x_k) \Longleftrightarrow
{\cal N} \models \varphi(x_1,\ldots,x_k).
\end{eqnarray} 

\begin{lemma}\label{folklore} \textbf{\em (Woodin)} Assume ${\sf MA}_{\omega_1}$, that ${\sf NS}_{\omega_1}$ is saturated, and that ${\cal P}(\omega_1)^\#$ exists. In the notation of Convention \ref{convention}:
\begin{enumerate}
\item[(1)] $g_A$ is a filter.
\item[(2)] If $g_A$ is ${\mathbb P}_{\rm max}$-generic over $L({\mathbb R})$,
then ${\cal P}({\omega_1}) \subset L({\mathbb R})[g]$.
\end{enumerate}
\end{lemma}

{\sc Proof.} This routinely follows from the proof of \cite[Lemma 3.12 and Corollary 3.13]{hugh} and from \cite[Lemmas 3.10 and 3.14]{hugh}. Let us sketch the argument. Let us first state the following.

\begin{claim}\label{claim_2.6}
Let $\theta \geq (2^{\aleph_1})^+$ be a cardinal, let $X \prec H_\theta$ be countable with $A \in X$, and let $\sigma \colon M \cong X$ be such that $M$ is transitive. Then 
\begin{enumerate}
\item[(a)] $\sigma^{-1}((H_{\omega_2};\in, {\sf NS}_{\omega_1} , A))$ is a ${\mathbb P}_{\rm max}$-condition.
\item[(b)] $\{ X \in {\cal P}(\omega_1^M) \cap M \colon \omega_1 \in \sigma(X) \}$ is $({\cal P}(\omega_1^M) \cap M) \setminus \sigma^{-1}({\sf NS})$-generic over $M$.
\item[(c)] For $i \leq \omega_1$ let $$X_i = {\rm Hull}^{H_{\theta}}(X \cup \sup \{ X_j \colon j<i \}) {\rm , }$$ let $\sigma_i \colon M_i \cong X_i$ be such that $M_i$ is transitive, and let, for $i \leq j \leq \omega_1$, $\pi_{ij} = \sigma^{-1}_j \circ \sigma_i$. Then $(M_i,\pi_{ij} \colon i \leq j \leq \omega_1)$ is a generic iteration of $$(M;\in,\sigma^{-1}({\sf NS}_{\omega_1}),A \cap \omega_1^M)$$
\item[(d)] $\sigma^{-1}((H_{\omega_2};\in, {\sf NS}_{\omega_1} , A)) \in g_A$.
\end{enumerate} 
\end{claim}

{\sc Proof} of Claim \ref{claim_2.6}. (a): This is by \cite[Lemmas 3.10 and 3.14]{hugh}. (b) and (c): This is by \cite[Lemma 3.12 and Corollary 3.13]{hugh}. (d): This follows immediately from (a) and (c). \hfill $\square$ (Claim \ref{claim_2.6})

\medskip
Let us now prove Lemma \ref{folklore}.

(1): Let ${\cal N}_0 = (N;\in,I,a) \in g_A$ as being witnessed by the generic iteration 
$({\cal N}_i , \sigma_{ij} \colon i \leq j \leq \omega_1)$. Let ${\cal M}_0 > {\cal N}_0$ as being witnessed by the generic iteration $({\cal M}_i , \pi_{ij} \colon i \leq j \leq \omega_1^{{\cal N}_0}) \in N$. Then $\sigma_{0 \omega_1}(({\cal M}_i , \pi_{ij} \colon i \leq j \leq \omega_1^{{\cal N}_0}))$ is easily seen to be a generic iteration of ${\cal M}_0$ which witnesses that ${\cal M}_0 \in g_A$.

Now let ${\cal N}_0^0 = (N^0;\in,I^0,a^0) \in g_A$ and ${\cal N}_0^1 = (N^1;\in,I^1,a^1) \in g_A$ as witnessed by the generic iterations 
${\cal I}^0 = ({\cal N}_i^0 , \sigma_{ij}^0 \colon i \leq j \leq \omega_1)$ and
${\cal I}^1 = ({\cal N}_i^1 , \sigma_{ij}^1 \colon i \leq j \leq \omega_1)$. Let $\sigma \colon M \cong X$ be as in Claim \ref{claim_2.6} with $\{ {\cal I}^0 , {\cal I}^1 \} \subset X$. Then ${\cal N}^0$ and ${\cal N}^1$ are both weaker than $\sigma^{-1}((H_{\omega_2};\in, {\sf NS}_{\omega_1} , A)) \in {\mathbb P}_{\rm max}$, cf.\ Claim \ref{claim_2.6} (a), as witnessed by the generic iterations $\sigma^{-1}({\cal I}^0) = ({\cal N}_i^0 , \sigma_{ij}^0 \colon i \leq j \leq \omega_1^M)$ and $\sigma^{-1}({\cal I}^1) = ({\cal N}_i^1 , \sigma_{ij}^1 \colon i \leq j \leq \omega_1^M)$.

(2): Let $Z \in {\cal P}(\omega_1)$. Let $\sigma \colon M \cong X$ be as in Claim \ref{claim_2.6} with $Z \in X$. Then $\sigma^{-1}((H_{\omega_2};\in, {\sf NS}_{\omega_1} , A)) \in g_A$ by Claim \ref{claim_2.6} (d)  and in fact if $(M_i,\pi_{ij} \colon i \leq j \leq \omega_1)$ is as in Claim \ref{claim_2.6} (c), then $Z \in M_{\omega_1}$, so that trivially $Z$ is also in the last iterate of $\sigma^{-1}((H_{\omega_2};\in, {\sf NS}_{\omega_1} , A))$ via the generic iteration which is the restriction of $(M_i,\pi_{ij} \colon i \leq j \leq \omega_1)$ to $\sigma^{-1}((H_{\omega_2};\in, {\sf NS}_{\omega_1} , A))$ and its images.
\hfill $\square$ (Lemma \ref{folklore})

\bigskip
%
Let $1 \leq k<\omega$, and let $D \in {\cal P}({\mathbb R}^k)$. We say that
$T$ is a tree on ${}^{k} \omega \times {\rm OR}$ iff $T \subset 
\bigcup_{n<\omega} ({}^n \omega)^k \times {}^n {\rm OR}$ and if $(s_0,...,s_{k-1},t) \in T$ and $m<\omega$, then $$(s_0 \upharpoonright m,...,s_{k-1} \upharpoonright m, t \upharpoonright m) \in T$$ We write $$[T] = \{ (x_0,...,x_{k-1},f) \colon \forall m<\omega \, (x_0 \upharpoonright m,...,x_{k-1} \upharpoonright m,f \upharpoonright m) \in T \}$$ and $p[T]$ for the projection of $T$, i.e., $$p[T] = \{ (x_0,...,x_{k-1}) \colon \exists f \, (x_0,...,x_{k-1},f) \in [T] \}$$

\begin{definition}
{\em The trees} $T$ {\em and} $U$ {\em on} ${}^{k} \omega \times {\rm OR}$
{\em witness that}
$D$ is {\em universally Baire} iff $D = p[T]$ and for all posets ${\mathbb P}$,
\begin{eqnarray}\label{complimenting_trees}
\Vdash_{\mathbb P} \, p[U] = {\mathbb R}^{k} \setminus p[T]{\rm . }
\end{eqnarray}
$D$ is called {\em universally Baire} iff there are trees $T$ and $U$ witnessing that 
$D$ is universally Baire. 

We denote by $\Gamma^\infty$ the collection of all $D \in \bigcup_{1 \leq k<\omega} \, {\cal P}({\mathbb R}^k)$ which are universally Baire.
\end{definition}

The concept of universally Baire set was isolated by Feng-Magidor-Woodin in \cite[Section 2]{FMW}; see also \cite[Definition 8.6]{book}. 

If $D \in \Gamma^\infty$, then there is an 
unambiguous version of $D$ in any forcing extension $V[g]$ of $V$, which as 
usual we 
denote by $D^*$ and which is equal to $p[T] \cap V[g]$ for some/all trees T and U which witness that $D$ is universally Baire. See \cite[p.\ 149f.]{book}.

We will call a pointclass consisting of universally Baire sets productive iff it is closed under complements and projections in a strong sense 
and 
for all $k<\omega$ and $D \in \Gamma \cap {\mathbb R}^{k+2}$,
\begin{eqnarray}\label{subtlety}
(\exists^{\mathbb R} D)^* = \{ {\vec x} \in {\mathbb R}^{k+1} \colon \exists y \in {\mathbb R} \, ({\vec x},y) \in D^* \}
\end{eqnarray}
will be true in every generic extension. The formal definition runs as follows.
 
\begin{definition}\label{defn-productive}
Let $\Gamma \subset \bigcup_{1 \leq k<\omega} \, {\cal P}({\mathbb R}^k)$. We say that $\Gamma$ is {\em productive} iff 
\begin{enumerate}
\item[(a)] $\Gamma \subset \Gamma^\infty$,
\item[(b)] for all $k<\omega$ and all $D \in \Gamma \cap {\cal P}({\mathbb R}^{k+1})$, ${\mathbb R}^{k+1} \setminus D \in \Gamma$ and if $k>0$, then
$\exists^{\mathbb R} D = \{ (x_0, \ldots , x_{k-1} ) \colon \exists x_k (x_0, \ldots , x_{k-1}, x_k) \in D \} \in \Gamma$, and
\item[(c)] for all $k<\omega$ and all $D \in \Gamma \cap {\cal P}({\mathbb R}^{k+2})$, if the trees $T$ and $U$ on ${}^{k+2} \omega \times {\rm OR}$ witness that $D$ is universally Baire 
and if 
\begin{eqnarray}\label{U-tilde}
{\tilde U} = \{ (s \upharpoonright (k+1), (s(k+1),t)) \colon (s,t) \in U \}{\rm , }
\end{eqnarray}
then 
there is a tree ${\tilde T}$ on ${}^{k+1} \omega \times {\rm OR}$ such that for all posets ${\mathbb P}$, 
\begin{eqnarray}\label{T-tilde}
\Vdash_{\mathbb P} \, p[{\tilde U}] = {\mathbb R}^{k+1} \setminus p[{\tilde T}] {\rm . }
\end{eqnarray}
\end{enumerate}
\end{definition}

While (c) of Definition \ref{defn-productive} canonically ensures that every productive pointclass is closed under projections, at least on the face of its definition, $\Gamma$
being productive is stronger than having that $\Gamma \subset \Gamma^\infty$ and
$\Gamma$ is closed under complements and projections (\cite[Question 3]{FMW} exactly asks if the former is really stronger than the latter). 

\begin{lemma}\label{absoluteness_for_productive}
If $\Gamma$ is productive and if $D \in \Gamma$, then any projective statement about $D$ is absolute between $V$ and any forcing extension of $V$,
i.e., if $\varphi$ is projective, $x_1$, $\ldots$, $x_k \in {\mathbb R}$, and ${\mathbb P}$ is any poset, then
$$V \models \varphi(x_1, \ldots , x_k, D) \ \Longleftrightarrow \ \Vdash_{\mathbb P} \, \varphi({\check x}_1, \ldots , {\check x}_k, {D}^*).$$ 
\end{lemma}

Lemma \ref{absoluteness_for_productive} is shown by a trivial induction on the complexity of $\varphi$. 

Let $e \colon {\mathbb R} \rightarrow {\rm HC}$ be a fixed simple coding of hereditarily countable sets by reals, see e.g.\ \cite[p.\ 179]{both}. A set $D \subset {\rm HC}$ is then called {\em universally Baire in the codes} iff the code set $\{ x \in {\mathbb R} \colon e(x) \in D \}$ of $D$ is universally Baire. If this is the case, then every forcing extension of $V$ will have its unique new version of $D$, which we denote by $D^*$. If the code set of $D$ is a member of a productive pointclass, then for every forcing ${\mathbb P}$,
\begin{eqnarray}\label{el_equivalent}
({\rm HC}; \in, D) \prec ({\rm HC}^{V^{\mathbb P}}; \in, D^*){\rm . }
\end{eqnarray}

A classical variant of Lemma \ref{absoluteness_for_productive} is Shoenfield's absoluteness theorem, see e.g.\ \cite[Corollary 7.21]{book}. It states that 
if $M \subset N$ are both transitive models of a sufficiently rich fragment of ${\sf ZFC}$ such that $\omega_1^V \subset M$, then 
\begin{eqnarray}\label{el_equivalent_shoenfield}
({\rm HC}^M; \in) \prec_{\utilde{\Sigma}_1} ({\rm HC}^{N}; \in){\rm , }
\end{eqnarray}
where (\ref{el_equivalent_shoenfield}) means that (\ref{ee}) holds true with $\varphi$ restricted to $\Sigma_1$ formulae (and ${\rm HC}^M$, ${\rm HC}^N$ playing the roles of ${\mathcal M}$, ${\mathcal N}$, respectively).

\cite[Question 3]{FMW}
is concerned
with the question about the connection of, on the one hand, projective absoluteness with respect to forcing extensions and, on the other hand, having that every projective set is universally Baire (see \cite[Questions 1 and 7]{FMW}).

\begin{theorem}\label{Gammainfty-productive} \textbf{\em  (Woodin)} Assume that there is a proper class of Woodin cardinals. Then $\Gamma^\infty$ is productive.
\end{theorem}

{\sc Proof.} A theorem of Woodin says that in the presence of a proper class of Woodin cardinals, every set in $\Gamma^\infty$ is weakly homogeneously Suslin,
see e.g.\ \cite[Theorem 3.3.8]{paul} and \cite[Theorem 1.2]{forcing-free}. 
Every tree ${\tilde U}$ witnessing that a given set $D$ of reals is weakly homogeneously Suslin comes with a canonical tree ${\tilde T}$ for 
${\mathbb R} \setminus D$ in such a way that ${\tilde U}$ and ${\tilde T}$ are connected as in (c) of Definition 
\ref{defn-productive}. For the construction of ${\tilde T}$ see e.g. \cite[p.\ 455]{kanamori}. \cite[Proposition 32.6]{kanamori} formulates how the two trees are connected. The main result of Martin and Steel
from \cite{proj-determinacy} is then that Woodin cardinals may be used to show that ${\tilde T}$ is homogeneous, cf.\ \cite[Theorem 32.11]{kanamori}.
That way, it follows that $\Gamma^\infty$ is productive provided that there is a proper class of Woodin cardinals. \hfill $\square$ (Theorem \ref{Gammainfty-productive})

\bigskip
For any set $X$, $M_\omega^\#(X)$ denotes the least active $X$-mouse which has infinitely many Woodin cardinals. See \cite{PWIM}.

\begin{theorem}\label{steel} \textbf{\em (Steel)} Assume {\sf PFA}. Then the universe is closed under the operation $X \mapsto M_\omega^\#(X)$. In particular, every set of reals in $L({\mathbb R})$ is universally Baire, and $\bigcup_{k<\omega} \, {\cal P}({\mathbb R}^k) \cap L({\mathbb R})$ is productive. 
\end{theorem}

{\sc Proof.} The proof from \cite{john} produces the result that under {\sf PFA}, the universe is closed under the operation $X \mapsto M_\omega^\#(X)$. 
The rest is given by standard inner model theoretic arguments, see e.g.\ \cite[Section 3, pp.\ 187f.]{both}. \hfill $\square$ (Theorem \ref{steel})

\bigskip
By Lemma \ref{folklore} and Theorem \ref{steel},
Theorem \ref{main_thm} 
follows from the following more general statement.

\begin{theorem}\label{main_thm_gen}
Let $\Gamma \subset \bigcup_{1 \leq k<\omega} \, {\cal P}({\mathbb R}^k)$. Assume that
\begin{enumerate}
\item[(i)] $\Gamma = \bigcup_{1 \leq k<\omega} \, {\cal P}({\mathbb R}^k) \cap L(\Gamma,{\mathbb R})$,
\item[(ii)] $\Gamma$ is productive, and
\item[(iii)] Martin's Maximum${}^{++}$ holds true.
\end{enumerate}
Then $g_A$ is ${\mathbb P}_{\rm max}$-generic over $L(\Gamma,{\mathbb R})$.
\end{theorem}

The abbreviation $(*)_\Gamma$ was introduced in \cite[Definition 4.1]{both} to denote a straightforward generalization of $(*)$ to larger pointclasses. For a pointclass $\Gamma \supset {\cal P}({\mathbb R}) \cap L({\mathbb R})$, $(*)_\Gamma$ is the statement that every set in $\Gamma$ is determined, and there is a filter $g \subset {\mathbb P}_{\rm max}$ which has nonempty intersection with every dense set (coded by a set) in $\Gamma$ and is such that ${\cal P}(\omega_1) \subset L({\mathbb R})[g]$.

\begin{corollary}
Assume that there is a proper class of Woodin cardinals.
Let $\Gamma \subset \bigcup_{1 \leq k<\omega} \, {\cal P}({\mathbb R}^k) \cap \Gamma^\infty$. Suppose that
(i)--(iii) from the statement of Theorem \ref{main_thm_gen} are satisfied.
Then $(*)_\Gamma$ holds true.
\end{corollary}


Theorem
\ref{main_thm_gen} readily follows from the following Lemma via a standard application
of ${\sf MM}^{++}$.

\begin{lemma}\label{main_claim}
Let $\Gamma \subset \bigcup_{1 \leq k<\omega} \, {\cal P}({\mathbb R}^k)$. Assume that
\begin{enumerate}
\item[(i)] $\Gamma = \bigcup_{1 \leq k<\omega} \, {\cal P}({\mathbb R}^k) \cap L(\Gamma,{\mathbb R})$,
\item[(ii)] $\Gamma$ is productive, and
\item[(iii)] ${\sf NS}_{\omega_1}$ is saturated.
\end{enumerate}

Let $D \subset {\mathbb P}_{\rm max}$ be open dense, $D \in L(\Gamma,{\mathbb R})$.
With $A$ being as in Convention \ref{convention},
there is then a stationary set preserving forcing ${\mathbb P}$ of size $2^{\aleph_2}$ such that in $V^{\mathbb P}$
there is some $p = {\cal N}_0 = (N;\in,I,a) \in D^*$ and some generic iteration $$({\cal N}_i,\sigma_{ij} \colon i \leq j \leq \omega_1)$$
of $p = {\cal N}_0$ of length $\omega_1+1$ such that if we write ${\cal N}_{\omega_1} = (N_{\omega_1};\in,I^*, a^*)$, then $I^* = ({\sf NS}_{\omega_1})^{V^{\mathbb P}} \cap N_{\omega_1}$ and $a^* = A$.
\end{lemma}

{\sc Proof} of Theorem \ref{main_thm_gen} from Lemma \ref{main_claim}. ${\sf MM}$ implies that ${\sf NS}_{\omega_1}$ is saturated, see \cite[Theorem 12]{FMS}. By Lemma \ref{folklore}, it remains to show that $D \cap g_A \not= \emptyset$ for every
open dense $D \subset {\mathbb P}_{\rm max}$, $D \in \Gamma$. Here, $g_A$ is as in Convention \ref{convention}.

Let us fix such $D$. The statement that there is a $p$ as in the conclusion of 
Lemma \ref{main_claim}, which is tantamount to saying that there is a $p \in D \cap g_A$, is easily seen to be $\Sigma_1$ expressible over the structure
$(H_{\omega_2};\in,{\sf NS}_{\omega_1},A,D)$. By the conclusion of Lemma \ref{main_claim}, the existence of such a $p$ may be forced by a stationary set preserving forcing.
Hence by ${\sf MM}^{++}$, cf.\ \cite[Theorem 10.124]{hugh}, 
there is such a $p$ in $V$. \hfill $\square$ (Theorem \ref{main_thm_gen})

\bigskip

As the proof of Theorem \ref{main_thm_gen} from Lemma \ref{main_claim} shows,
we don't need the full power of ${\sf MM}^{++}$ 
in order to derive Theorem \ref{main_thm_gen} from Lemma \ref{main_claim}.
Instead, a bounded version of ${\sf MM}^{++}$ suffices; it may be defined as follows, see
\cite[Definition 10.123]{hugh}. 

\begin{definition}
For $D \in \Gamma^\infty$, $D$-${\sf BMM}^{++}$ is the statement that for all ${\mathbb P}$ which are stationary set preserving, $$(H_{\omega_2}^V;\in,{\sf NS}_{\omega_1}^V,D) \prec_{\Sigma_1} (H_{\omega_2}^{V^{\mathbb P}};\in,{\sf NS}_{\omega_1}^{V^{\mathbb P}},D^*).$$
\end{definition}

Modulo large cardinals, $(*)$ is then actually {\em equivalent} to $D$-${\sf BMM}^{++}$ for all $D \in {\mathcal P}({\mathbb R}) \cap L({\mathbb R})$. Let us first
state a more general fact, Theorem \ref{2.9}, which gives the characterization of $(*)$, i.e.\ Theorem \ref{char_star}, as a special case. 

\begin{theorem}\label{2.9}
Assume that there is a proper class of Woodin cardinals.
Let $\Gamma \subset \bigcup_{1 \leq k<\omega} \, {\cal P}({\mathbb R}^k)$. Assume that
\begin{enumerate} 
\item[(i)] $\Gamma = \bigcup_{1 \leq k<\omega} \, {\cal P}({\mathbb R}^k) \cap L(\Gamma,{\mathbb R})$,
\item[(ii)] $\Gamma$ is productive.
\end{enumerate}
The following statements are then equivalent, with $g_A$ being as in Convention \ref{convention}.
\begin{enumerate}
\item[(1)] $D$-${\sf BMM}^{++}$ holds true for all $D \in \Gamma$. 
\item[(2)] $g_A$ is ${\mathbb P}_{\rm max}$-generic over $L(\Gamma,{\mathbb R})$.
\end{enumerate}
\end{theorem}

{\sc Proof.} (2) $\Longrightarrow$ (1): This is exactly by the proof of (A) $\Longrightarrow$ (B) of \cite[Theorem 2.7]{david-ralf}. 

(1) $\Longrightarrow$ (2): We may first force ${\sf NS}_{\omega_1}$ to be saturated
by a stationary set preserving forcing, see e.g.\ \cite[Theorem 2.64]{hugh}.
The rest is then by the proof of Theorem \ref{main_thm_gen} from Lemma \ref{main_claim}
which was given above.
 \hfill $\square$ (Theorem \ref{2.9})
 

\begin{theorem}\label{char_star} Assume that there is a proper class of Woodin cardinals. 
The following statements are then equivalent.
\begin{enumerate}
\item[(1)] $D$-${\sf BMM}^{++}$ holds true for all 
$D \in {\cal P}({\mathbb R}) \cap L({\mathbb R})$.
\item[(2)] $(*)$.
\end{enumerate}
\end{theorem}

Let us now give a {\sc Proof} of Theorem \ref{char_of_star} from Lemma \ref{main_claim}.
(1) $\Longrightarrow$ (2) is weaker than (1) $\Longrightarrow$ (2) of Theorem
\ref{Wo1}.
Let us now assume (2) and show (1). Fix $D \subset {\mathbb P}_{\rm max}$, any open dense set in $L({\mathbb R})$. As the statement of the theorem assumes a supercompact cardinal to exist, there is a semi-proper (and hence stationary set preserving) forcing ${\mathbb P}$ such that ${\sf MM}^{++}$ holds true in $V^{\mathbb P}$. 
Inside $V^{\mathbb P}$, we will have that ${\sf MM}^{++}$ yields via
Lemma \ref{main_claim} that for all $A' \subset \omega_1$ with $\omega_1^{L[A']} =
\omega_1$ there will be some $p = (N;\in,I^*,a^*) \in D^*$ and some generic iteration $$({\cal N}_i,\sigma_{ij} \colon i \leq j \leq \omega_1)$$
of $p = {\cal N}_0$ of length $\omega_1+1$ such that if we write ${\cal N}_{\omega_1} = (N_{\omega_1};\in,I^*, a^*)$, then $I^* = ({\sf NS}_{\omega_1})^{V^{\mathbb P}} \cap N_{\omega_1}$ and $a^* = A'$. This is a statement which is $\Pi_2$ over the
structure mentioned in (2) of Theorem \ref{char_of_star}.
This statement will therefore be true in $V$, which readily implies that
$g_A$ is ${\mathbb P}_{\rm max}$-generic over $L({\mathbb R})$
and ${\cal P}(\omega_1) \subset L({\mathbb R})[g_A]$, where $g_A$ is as in Convention \ref{convention}. \hfill $\square$ (Theorem \ref{char_of_star})

\bigskip
The forcing which we designed in order to produce Lemma \ref{main_claim} is a souped up version of the forcings from \cite{cs} and \cite{ds}, which are in turn variants of the ${\cal L}$-forcing of Jensen as being developed e.g.\ in \cite{ronald_l_forcing}.\footnote{One of the referees informs us that J.\ Keisler in \cite{kei71} and \cite{kei73} 
developed forcings which work in a similar fashion.}  
All these forcings may be construed as building uncountable models as term models of a given language, ${\cal L}$, with the forcing conditions being finite fragments of a consistent and complete ${\cal L}$-theory which will give those term models, augmented by ``side conditions'' which will guarantee that the forcing only collapses cardinals in a controlled way. Our forcing will change the cofinalities of $\omega_2$ and $\omega_3$ to $\omega$ and $\omega_1$, respectively, and it won't collapse any other cardinal
outside of the (possibly empty) half-open interval $(\omega_3,2^{\aleph_2}]$. 

Let us give an outline of the proof of Lemma \ref{main_claim}.

To prove Lemma \ref{main_claim}, we aim to build a stationary set preserving forcing $\mathbb P$ which adds a generic iteration of some ${\mathbb P}_{\rm max}$-condition $(N; \in, I, a)$ coded by a real in the projection of a tree ${\tilde T}$ projecting to the set of codes for conditions in our given dense set $D$. Moreover, we want this iteration to send the distinguished set $a$ of $(N; \in, I, a)$ to $A$, and we want every $I^*$-positive set in the final model $(N^*;\in,I^*,A)$ to be a stationary subset of $\omega_1$ in $V^{\mathbb P}$. Our approach is to think of all the relevant objects -- $(N; \in, I, a)$, a branch through ${\tilde T}$ projecting to a real coding $(N; \in, I, a)$, and the generic iteration of $(N; \in, I, a)$ of length $\omega_1+1$ -- as being given by ``term models'' in a suitable language, ${\cal L}$, and add them via finite approximations.
Thus, the working parts of our forcing will be finite sets $p$ of sentences from $\mathcal L$ providing partial information about the above objects. We will require these finite bits of information $p$ to be realized in some outer model.\footnote{$W$ is an outer model iff $W$ is a transitive model of {\sf ZFC} with
$W \supset V$ and which has the same ordinals as $V$; in other words, $W$ is an 
outer model iff $V$ is an inner model of $W$.} 
The existence of such an outer model will be absolute to
any generic extension of $V$ via ${\rm Col}(\omega, \omega_2)$.

In 
$V^{{\rm Col}(\omega, \omega_2)}$, $$(H_{\omega_2}^V; \in, {\sf NS}_{\omega_1}^V,A)$$ becomes a ${\mathbb P}_{\rm max}$-condition and $p[{\tilde T}]=D^*$ is still dense, so that in $V^{{\rm Col}(\omega, \omega_2)}$
there is a ${\mathbb P}_{\rm max}$-condition $(N; \in, I, a) \in D^*$ which is 
stronger than $(H_{\omega_2}^V; \in, {\sf NS}_{\omega_1}^V,A)$.
We may now iterate $(N; \in, I, a)$ in length $\omega_1^{V^{{\rm Col}(\omega, \omega_2)}}+1=\omega_3^V+1$ so as to produce $$\sigma \colon (N;\in,I,a) \rightarrow
(N^*;\in,I^*,a^*).$$ 
If $(M_i,\pi_{ij} \colon i \leq j \leq \omega_1^N) \in N$ is the generic iteration of $(H_{\omega_2}^V; \in, {\sf NS}_{\omega_1}^V,A)$ witnessing that $(N; \in, I, a)$ is stronger than 
$(H_{\omega_2}^V; \in, {\sf NS}_{\omega_1}^V,A)$, then
$\sigma((M_i,\pi_{ij} \colon i \leq j \leq \omega_1^N)) = (M_i,\pi_{ij} \colon i \leq j \leq 
\omega_3^V)$ is an extension of that iteration. We have that 
$$M_{0} = (H_{\omega_2}^V; \in, {\sf NS}_{\omega_1}^V,A){\rm , }$$ and
$$\pi_{0, \omega_3^V}
\colon M_{0} \rightarrow
M_{\omega_3^V}$$ may be
lifted to a generic iteration $${\tilde \pi} \colon V\to M$$ of $V$, for a transitive $M$, such that ${\tilde \pi} \supset \pi_{0, \omega_3^V}$ and ${\tilde \pi}(M_0)=M_{\omega_3^V}$. See \cite[Lemma 3.8]{hugh}.

\bigskip
\begin{tikzpicture}

\def\x{2.1};

\def\y{0.5};

\node (pT) at (0*\x,  0*\y){$p[{\tilde T}]$};

\node (ppiT) at (0*\x+1.5, 0*\y){$p[{\tilde \pi}({\tilde T})]$};

\node (N) at (0*\x, -2*\y) {$(N;\in, I, a)$};

\node (Nast) at (2*\x, -2*\y) {$(N^*;\in, I^*, a^*)$};

\node (M0) at (-2*\x, -4*\y) {$M_0$};

\node (MN) at (0*\x, -4*\y) {$M_{\omega_1^N}$};

\node (M3) at (2*\x, -4*\y) {$M_{\omega_3^V}$};

\node (H) at (-2*\x, -6*\y) {$(H_{\omega_2}^V;\in, \mathsf{NS}_{\omega_1}^V, A)$};

\node (V) at (-2*\x, -8*\y) {$V$};

\node (Mp) at (2*\x, -8*\y) {$M$};


\path (pT)--(ppiT) node[midway]{$\subseteq$};

\path (N)--(pT) node[midway, rotate=90]{$\in$};

\draw[->] (N)--(Nast) node[midway, above]{$\sigma$};

\path (MN)--(N) node[midway, rotate=90]{$\in$};

\path (M3)--(Nast) node[midway, rotate=90]{$\in$};

\draw[->] (M0)--(MN) node[midway, above]{$\pi_{0 \omega_1^N}$};

\draw[->] (MN)--(M3) node[midway, above]{$\pi_{\omega_1^N \omega_3^V}$};

\path (H)--(M0) node[midway, rotate=90]{$=$};

\path (M3)--(Mp) node[midway, rotate=270]{$\in$};

\path (H)--(V) node[midway, rotate=270]{$\in$};

\draw[->] (V)--(Mp) node[midway, above]{${\tilde \pi}$};


\end{tikzpicture}

\bigskip
We now see that $V^{{\rm Col}(\omega, \omega_2)}$ contains objects like the ones we intend to add by our forcing---namely $(N; \in, I, a)$, a branch through ${\tilde T}$ projecting to a real coding $(N; \in, I, a)$, and the generic iteration of $(N; \in, I, a)$ of length $\omega_1+1$---albeit not defined relative to the parameters ${\tilde T}$, $H_{\omega_2}^V$, and $A$, but relative to ${\tilde \pi}({\tilde T})$, $H_{\omega_2}^M$, and ${\tilde \pi}(A)$. 
The statement that such objects exist is $\Sigma_1$ in the parameters $H_{\omega_2}^M$, ${\tilde \pi}(A)$, and a Skolem hull\footnote{See Claim \ref{ex_of_T}.} of ${\tilde \pi}({\tilde T})$ of size $\aleph_2^M$, which will both be elements of ${\rm HC}^{M^{{\rm Col}(\omega,{\tilde \pi}(\omega_2))}}$. 
By Shoenfield absoluteness (\ref{el_equivalent_shoenfield}), see e.g.\ \cite[Corollary 7.21]{book}, 
such objects will also exist in $M^{{\rm Col}(\omega, {\tilde \pi}(\omega_2))}$. 

The point that $\pi_{0,\omega_3^V}$ could be lifted to ${\tilde \pi}$ is then the following.
The statement that objects like the ones we intend to add by our forcing exist in
$M^{{\rm Col}(\omega, {\tilde \pi}(\omega_2))}$
may now be pulled back via ${\tilde \pi}$. This buys us that objects like the ones we intend to add by our forcing exist in $V^{{\rm Col}(\omega, \omega_2)}$ -- and this time with the right parameters ${\tilde T}$, $H_{\omega_2}^V$, and $A$. The argument that combined lifting $\pi_{0,\omega_3^V}$ to ${\tilde \pi}$, applying Shoenfield absoluteness, and pulling back the statement of interest was crucial to arrive at the desired conclusion, viz.\ that objects like the ones we intend to add by our forcing exist in $V^{{\rm Col}(\omega, \omega_2)}$.\label{discussion_why_M} This will be our starting point for cooking up the forcing ${\mathbb P}$.


In order to prove that our forcing ${\mathbb P}$ preserves stationary subsets of $\omega_1$ we will need an argument exploiting lifting, Shoenfield absoluteness, and pulling back. In order to be able to run this argument we will need our forcing to approximate, not only  the objects we are ultimately interested in obtaining, but also the iteration
$(M_i,\pi_{ij} \colon i \leq j \leq \omega_1^N) \in N$. (See footnote \ref{footnote}.) 
We will think of the objects themselves, which exist in $V^{{\rm Col}(\omega,\omega_2)}$, as ``certificates''  for some finite piece of information about them. The idea is then to have our forcing consist of finite sets of $\mathcal L$-sentences for which there is a ``certificate'' in $V^{{\rm Col}(\omega,\omega_2)}$.

The problem with the above strategy is that, although a forcing $\mathbb P$ like the one we have described would in fact add the desired objects, one would still need to show that it preserves stationary subsets of $\omega_1$ and that every positive set in the final model of the iteration being added by $\mathbb P$ is in fact stationary in that extension. 
Our forcing ${\mathbb P}$ will be a subset of $H_{\omega_3}$, and 
one tool for taking care of these issues is the use of a diamond sequence $\langle (Q_\lambda, A_\lambda) \,:\,\lambda < \omega_3\rangle$ consisting of transitive structures\footnote{I.e., the underlying universe $Q_\lambda$ will be transitive and $A_\lambda \subset Q_\lambda$ will be a distinguished predicate of the structure.} in $H_{\omega_3}$ in order to guess $(H_{\omega_3},{\dot C})$, where ${\dot C}$ is a $\mathbb P$-name for a 
club in $\omega_1$, ${\dot C} \subset H_{\omega_3}$. 
That $(H_{\omega_3},{\dot C})$ be guessed means that there are stationarily many $\lambda<\omega_3$ such that $(Q_\lambda,A_\lambda)$ is an elementary substructure of $(H_{\omega_3},{\dot C})$. See $(\Diamond)$ on p.\ \pageref{defn_diamond}.

Imagine that ${\mathbb P}$ is a forcing which adds the desired objects, but which also preserves stationary subsets of $\omega_1$. Let ${\dot C}$ be a $\mathbb P$-name for a club in $\omega_1$,
${\dot C} \subset H_{\omega_3}$, and let $S \subset \omega_1$ be stationary in $V$.
Let $g$ be ${\mathbb P}$-generic 
over $V$. There will be some $\lambda<\omega_3$ such that $(H_{\omega_3},{\dot C})$ is guessed by $(Q_\lambda,A_\lambda)$ and in $V[g]$ there will be some countable elementary substructure
$X$ of $(Q_\lambda,A_\lambda)$ such that, writing $X[g] = \{ \tau^g \colon \tau \in V^{\mathbb P} \cap X \}$,
\begin{enumerate}
\item[(a)] $X \cap \omega_1 \in S$, and
\item[(b)] $X[g] \cap Q_\lambda = X \cap Q_\lambda$.
\end{enumerate}
That $(Q_\lambda,A_\lambda)$ is an elementary substructure of $(H_{\omega_3},{\dot C})$ will then mean in practice that $X \cap \omega_1 \in {\dot C}^g$, so that
$X \cap \omega_1$ witnesses that $S$ is still stationary in $V[g]$.
Calling some $g$ with the property (b) ``$({\mathbb P},X)$-generic''
there is, however, no reason to expect an $X$ such that $g$ is $({\mathbb P},X)$-generic to exist in $V$ (in fact, $V$ won't have such $X$). 

When defining ${\mathbb P}$, we will turn this around and have our conditions also approximate finite bits of information about such elementary substructures $X$. 

Our key tool for taking care of the above issue is then to define ${\mathbb P}$ as the last forcing from a
recursively defined
$\subset$-increasing sequence $\vec{\mathbb P} = (\mathbb P_\lambda\,:\,\lambda \leq \omega_3
)$. Each $\mathbb P_\lambda$ will be a subset of $Q_\lambda$. 
Hence, when defining ${\mathbb P}_{\eta}$, $\eta \leq \omega_3$,
if $\lambda<\eta$, then we already know what it means for some $g$ to be (partially) ${\mathbb P}_\lambda$-generic over $(Q_\lambda, A_\lambda)$, and if
$X$ is a countable elementary substructure of $(Q_\lambda, A_\lambda)$, then $X[g]$ may be assigned a meaningful interpretation as
$X[g] = \{ \tau^g \colon \tau \in V^{{\mathbb P}_\delta} \cap X \}$.
We will maintain that 
at each stage $\eta$ in the construction of $\vec{\mathbb P}$ we define $\mathbb P_\eta$ by saying that a finite set $p$ of $\mathcal L\cap Q_\eta$-sentences is in $\mathbb P_\eta$ if and only if there is a certificate for $p$ which extends $p$ and when intersected with each of the side conditions $X_{\lambda} \prec (Q_\lambda , A_\lambda)$ (also given by the certificate), for $\lambda<\eta$, is generic over $X_\eta$ for the already defined forcing $\mathbb P_\lambda$. 
The role of condition (b) above (with $X_\lambda$ replacing $X$) will be that if $A_\lambda$ codes the name of a club subset ${\dot C}$ of $\omega_1$, then some $p \in g$ must force that $X_\lambda \cap \omega_1 \in {\dot C}$, so that in the light of (a) above, $p$ also forces that $S$ has non-empty intersection with ${\dot C}$.

Our next section is entirely devoted to a proof of Lemma \ref{main_claim}.

\section{The forcing.}

Recall our Convention \ref{convention} which we are now going to make use of without further notice.
Let us assume throughout the hypotheses of Lemma \ref{main_claim}.
We aim to verify its conclusion.

Let us fix $D \subset {\mathbb P}_{\rm max}$, an open dense set in $L(\Gamma,{\mathbb R})$. The fact that $D$ is open dense may be written as $$\forall p \in {\mathbb P}_{\rm max} \, \exists q \leq p, \, q \in D \wedge \forall p \in D \, \forall q \leq p, \, q \in D.$$ 
By hypothesis (ii) in the statement of Lemma \ref{main_claim} (i.e., $\Gamma \subset \Gamma^\infty$ is productive), we 
may apply Lemma \ref{absoluteness_for_productive} to conclude that 
 (\ref{el_equivalent}) on p.\ \pageref{el_equivalent} holds true with $D$ and $\mathbb P={\rm Col}(\omega,\omega_2)$, i.e., $$({\rm HC} ; \in, D) \prec ({\rm HC}^{V^{{\rm Col}(\omega,\omega_2)}} ; \in, D^*).$$ This will ensure that
\begin{enumerate}
\item[(D.1)] $V^{{\rm Col}(\omega,\omega_2)} \models$ ``$D^*$ is an open dense
subset of ${\mathbb P}_{\rm max}$.''
\end{enumerate}
Let us identify $D$ with a canonical set of reals coding the elements of $D$,\footnote{We will later have to spell out a bit more precisely in which way we aim to have the elements of $p[T]$ code the elements of $D$, see (C.2) and ($\Sigma$.5) below.} and
let ${\tilde T} \in V$ be a
tree on $\omega \times \theta$, for some ordinal $\theta$, such that 
\begin{enumerate}
\item[(D.2)] $V^{{\rm Col}(\omega,\omega_2)} \models D^*=p[{\tilde T}]$.\footnote{An easy Skolem hull argument may be used to show that we might actually pick ${\tilde T} \in V$ as a tree
on $\omega \times 2^{{\aleph_2}}$. We won't need that, though, but we shall prove and make use of a related fact below, see (\ref{q-star}).}
\end{enumerate}

Let
$h$ be ${\rm Col}(\omega,\omega_2)$-generic over $V$. Inside $V[h]$, 
\begin{eqnarray}\label{its-a-Pmax-condition!}
((H_{\omega_2})^V;\in, ({\sf NS}_{\omega_1})^V , A)
\end{eqnarray} is a ${\mathbb P}_{\rm max}$ condition, call it $p$.
Let $q^* \in ({\mathbb P}_{\rm max})^{V[h]}$, $q^* < p$, $q^* \in D^*$, cf.\ (D.1).
Let $q^* = (N^*; \in , I^* , a^*)$. Identifying $q^*$ with some real coding it, we have that
$q^* \in p[{\tilde T}]$, cf.\ (D.2). 

\begin{claim}\label{ex_of_T}
There is a tree ${T} \in V$ on $\omega \times \omega_2$ such that 
\begin{eqnarray}\label{q-star}
q^* \in p[{T}] \subset p[{\tilde T}].
\end{eqnarray} 
\end{claim}
{\sc Proof} of Claim \ref{ex_of_T}. 
Let $q^* = \sigma^h$, where $\sigma \in V^{{\rm Col}(\omega,\omega_2)}$. We may assume that $\sigma \in H_{\omega_3}$. Recall that ${\tilde T}$ is on $\omega \times \theta$. Let $X \in V$, $X \prec H_{\theta^+}^V$ be such that $\omega_2+1 \cup \{ \sigma , {\tilde T} \} \subset X$ and ${\rm Card}(X)=\aleph_2$. Let $\pi \colon P \cong X \prec H_{\theta^+}^V$ be such that $P$ is transitive, and write $T = \pi^{-1}({\tilde T})$. We have that $\pi(\sigma)=\sigma$, and $\pi$ lifts to ${\tilde \pi} \colon P[h] \rightarrow 
H_{\theta^+}^{V[h]}$ with ${\tilde \pi}(q^*) = {\tilde \pi}(\sigma^h) = \pi(\sigma)^h = \sigma^h = q^*$. As $q^* \in p[{\tilde T}]$, the elementarity of ${\tilde \pi}$ then yields that $q^* \in p[T]$. The tree $T$ is on $\omega \times P \cap {\rm OR}$, but using a bijection of $P \cap {\rm OR}$ with $\omega_2$, we may construe it as a tree on $\omega \times \omega_2$. 
\hfill $\square$ (Claim \ref{ex_of_T})

\bigskip
Let us fix $T$ as in Claim \ref{ex_of_T}.
Let us write 
\begin{eqnarray} \kappa = \aleph_3 {\rm , }
\end{eqnarray}
so that $T \in H_\kappa$.
Let $d$ be ${\rm Col}(\kappa,\kappa)$-generic over $V$. In $V[d]$, let
$({\bar A}_\lambda \colon \lambda<\kappa)$ be a $\Diamond_\kappa$-sequence, i.e., 
for all ${\bar A} \subset \kappa$, $\{ \lambda<\kappa \colon {\bar A} \cap \lambda = {\bar A}_\lambda \}$ is stationary.
Also, let $c \colon \kappa \rightarrow H_\kappa^{V} = H_\kappa^{V[d]}$,
$c \in V[d]$, be bijective.
For $\lambda < \kappa$, let 
\begin{eqnarray}
Q_\lambda = c \mbox{''} \lambda \ \mbox{  and  } \ A_\lambda = c \mbox{''} {\bar A}_\lambda.
\end{eqnarray} 
An easy closure argument will give us some club\label{defn-club} $C \subset \kappa$ such that for all $\lambda \in C$,
\begin{enumerate}
\item[(i)] $Q_\lambda$
is transitive, 
\item[(ii)] $\{ T , ((H_{\omega_2})^V ; \in , ({\sf NS}_{\omega_1})^V ,A) \} \cup
(\omega_2+1) \subset Q_\lambda$,
\item[(iii)] $Q_\lambda \cap {\rm OR} = \lambda$ (so that $c \upharpoonright \lambda \colon
\lambda \rightarrow Q_\lambda$ is bijective), and 
\item[(iv)] $(Q_\lambda;\in) \prec
(H_\kappa;\in)$.
\end{enumerate}
(ii) is true for all sufficiently large $\lambda < \kappa$, and (iv) is true for all $\lambda$ such that $Q_\lambda = c \mbox{''} \lambda$ is closed under some fixed set of Skolem functions for $H_\kappa$. As the set of $\lambda < \kappa$ with (i) and (iii) is each easily seen to be club, a club of $\lambda$ with the above properties certainly exists.
We will fix from now on some club $C \subset \kappa$ with (i) through (iv) for all $\lambda \in C$.

In $V[d]$, for all $P$, $B \subset H_\kappa$, the set of all $\lambda \in C$ such that
$$(Q_\lambda;\in,P \cap Q_\lambda , B \cap Q_\lambda) \prec (H_\kappa;\in,P,B)$$
is club, and the set of all $\lambda \in C$ such that $B \cap Q_\lambda = A_\lambda$ is stationary, so that 
\begin{enumerate}
\item[($\Diamond$)] For all $P$, $B \subset H_\kappa$ the set \label{defn_diamond}
$$\{ \lambda \in C \colon (Q_\lambda;\in , P \cap Q_\lambda , A_\lambda) \prec (H_\kappa;\in , P , B) \}$$
is stationary.
\end{enumerate}

We shall sometimes also write $Q_\kappa=H_\kappa$.
Readers who are familiar with Jensen's diamond will 
easily see that the principle which we refer to as ($\Diamond$) is actually equivalent to $\Diamond_\kappa$; see e.g.\ \cite[Definition 5.34]{book}. We shall use ($\Diamond$) to guess information about names for club subsets of $\omega_1$; 
this will play a crucial role in the verification that our forcing preserves stationary subsets of $\omega_1$.

\subsection{The definition of the forcing.}\label{3.1}

We shall now go ahead and 
produce a stationary set preserving forcing ${\mathbb P} \in V[d]$ of size $\kappa$
which adds some
$p \in D^*$ and some generic iteration $$(N_i,\sigma_{ij} \colon i \leq j \leq \omega_1)$$
of $p = N_0$ such that if we write $N_{\omega_1} = (N_{\omega_1};\in,I^*, a^*)$, then $I^* = ({\sf NS}_{\omega_1})^{V[d]^{\mathbb P}} \cap N_{\omega_1}$ and $a^* = A$. As the forcing ${\rm Col}(\kappa,\kappa)$ which added $d$ is certainly stationary set preserving, this will verify Lemma \ref{main_claim}.

${\sf NS}_{\omega_1}$ is still saturated in $V[d]$. This is true simply because 
forcing with ${\rm Col}(\kappa,\kappa)$ doesn't add any sequences of elements of $V$ of length $\aleph_2$. Moreover, (D.1) and (D.2) are true with $V$ being replaced by $V[d]$, as no reals are added and ${\mathbb P}_{\rm max}$ remains unchanged. Hence, in order to simplify our notation, we shall in what follows write $V$ for $V[d]$, i.e., assume that, in addition to ``${\sf NS}_{\omega_1}$ is saturated'' 
plus (D.1) and (D.2), ($\Diamond$) is also true in $V$. 

Working under these hypotheses, we shall now
recursively define a $\subset$-increasing and continuous chain of forcings ${\mathbb P}_\lambda$ for all 
$\lambda \in C \cup \{ \kappa \}$.
 The forcing ${\mathbb P}$ will be ${\mathbb P}_\kappa$. The conditions in each ${\mathbb P}_\lambda$ will be finite sets of formulae of an associated first order language, ${\cal L}^\lambda$, which will be defined below. The order of each ${\mathbb P}_\lambda$ will be just reverse inclusion, i.e.,
$q \leq_{{\mathbb P}_\lambda} p$ iff $q \supset p$ for $p$, $q \in {\mathbb P}_\lambda$.\footnote{Every ${\mathbb P}_{\lambda}$ will be designed to add certain objects by means of finite sets of formulae giving partial information on these objects. As indicated by the description at the end of Section \ref{section2}, our intended objects are mentioned by the first forcing ${\mathbb P}_{\min(C)}$; the additional objects mentioned by the latter forcings are introduced in order to ensure that the final member of the sequence, i.e., ${\mathbb P}_\kappa$,  
has the desired property of preserving stationary sets and adding a correct generic iteration. (So it is natural to think of ${\mathbb P}_{\min(C)}$ as supporting the ``working part'' of the conditions in ${\mathbb P}_\kappa$, and all latter forcings as supporting also ``side conditions.'')  It might therefore be useful, on a first reading of the remainder of this section, to only focus on the definition and analysis of ${\mathbb P}_{\min(C)}$---in other words, to ignore everything involving clauses (C.6) to (C.8) in the definition below of potential certificate, to also ignore 
clauses ($\Sigma$.6) to ($\Sigma$.8) in the definition below of potential certificates being pre-certified,
and to skip Subsection \ref{preservation-subsection}---and to pay attention to the entirety of the present section only on a second reading.}

Assume that $\lambda \in C \cup \{ \kappa \}$ and ${\mathbb P}_\mu$ has already 
been defined in such a way that ${\mathbb P}_\mu \subset Q_\mu$
for all $\mu \in C \cap \lambda$. We aim to define ${\mathbb P}_\lambda$.

\begin{convention}\label{convention_real_code}
We say that $x \subset \omega$ is a real code for $N_0 = (N_0;\in,I,a)$
if
there is some surjection $f \colon \omega \rightarrow N_0$ such that
$x$ is 
the monotone enumeration 
of the G\"odel numbers of all expressions of the form $\ulcorner {\dot N}_0 \models \varphi({\dot n}_1, \ldots , {\dot n}_\ell, {\dot a}, {\dot I})
\urcorner$ such that 
$\varphi$ is a first order formula of the language associated to $(N_0;\in,I,a)$\footnote{I.e., the language of set theory augmented by predicates for $I$ and $a$}
and 
${N}_0 \models \varphi(f(n_1), \ldots , f(n_\ell),a,I)$ holds true. 
\end{convention}

We shall be interested in objects ${\mathfrak C}$ which exist in some outer model 
and which have the following properties. Any ${\mathfrak C}$ will be a triple of sets which are indexed by $\omega_1$, $\omega$, and $K$, respectively, with $K$ being a subset of $\omega_1$:
\begin{eqnarray}\label{certificate}
{\mathfrak C} = \langle \langle M_i , \pi_{ij} , N_i , \sigma_{ij} \colon i \leq j \leq \omega_1 \rangle , 
\langle (k_n , \alpha_n ) \colon n<\omega \rangle ,
\langle \lambda_\delta , X_\delta \colon \delta \in K \rangle \rangle {\rm , }
\end{eqnarray}
where
\begin{enumerate}
\item[(C.1)] $M_0$, $N_0 \in {\mathbb P}_{\rm max}$,
\item[(C.2)] $x = \langle k_n \colon n<\omega \rangle$ is a real code for $N_0 = (N_0;\in,I,a)$ in the sense of Convention \ref{convention_real_code} and
$x$ is 
$\langle (k_n , \alpha_n ) \colon n<\omega \rangle \in [T]$,
\item[(C.3)] $\langle M_i , \pi_{ij} \colon i \leq j \leq \omega_1^{N_0} \rangle \in N_0$ is a generic iteration of $M_0$ which witnesses that $N_0
< M_0$ in ${\mathbb P}_{\rm max}$,
\item[(C.4)] $\langle N_i , \sigma_{ij} \colon i \leq j \leq \omega_1 \rangle$ is a generic iteration of $N_0$ such that if $$N_{\omega_1} = (N_{\omega_1};\in,I^*,A^*){\rm , }$$ then 
$A^* = A$,\footnote{There is no requirement on $I^*$ matching the non-stationary ideal of some model in which ${\mathfrak C}$ exists.\label{no-requirement}}
\item[(C.5)] $\langle M_i , \pi_{ij} \colon i \leq j \leq \omega_1 \rangle = \sigma_{0 \omega_1}( \langle M_i , \pi_{ij} \colon i \leq j \leq \omega_1^{N_0}\rangle )$ and \begin{eqnarray}\label{Momega1}
M_{\omega_1} = ((H_{\omega_2})^V ; \in , ({\sf NS}_{\omega_1})^V ,A){\rm , }\footnotemark
\end{eqnarray}
\footnotetext{In particular, the distingushed ideal of $M_{\omega_1}$ is the true nonstationary ideal of $V$.}
\item[(C.6)] $K \subset {\omega_1}$, 
\end{enumerate}
and for all $\delta \in K$,
\begin{enumerate}
\item[(C.7)] $\lambda_\delta \in \lambda \cap C$,\footnote{Recall that $\lambda$ is the index of the forcing ${\mathbb P}_\lambda$ which we are about to define.} and if $\gamma<\delta$ is in $K$, then 
$\lambda_\gamma<\lambda_\delta$ and 
$X_\gamma \cup \{ \lambda_\gamma \} \subset X_\delta$, and
\item[(C.8)] $X_\delta \prec (Q_{\lambda_\delta};\in,{\mathbb P}_{\lambda_\delta},A_{\lambda_\delta})$ and $X_\delta \cap \omega_1 = \delta$.
\end{enumerate}

For future purposes, let us refer to any object ${\mathfrak C}$ as in (\ref{certificate}) which satisfies the above properties (C.1) through (C.8) as a {\em potential certificate}. As this definition of potential certificates 
has $A$, $T$, $\langle A_\mu\,:\,\mu\in C \cap \lambda\rangle$, and $\langle \mathbb P_\mu\,:\,\mu\in C\cap\lambda\rangle$ as hidden parameters, we shall also refer to them as {\em potential certificates relative to} $A$, $T$, $H_{\omega_2}^V$, $\langle A_\mu\,:\,\mu\in C \cap \lambda\rangle$, and $\langle \mathbb P_\mu\,:\,\mu\in C\cap \lambda \rangle$. 

We need to define a first order language ${\cal L}$ (independently from $\lambda$)
whose formulae will be able to describe some ${\mathfrak C}$ with the above properties by producing the 
models $M_i$ and $N_i$, $i <\omega_1$, as term models out
of equivalence classes of terms of the form ${\dot n}$, $n<\omega$.
The language ${\cal L}$ will have the following constants.
There will be one constant for each set in $H_\kappa$; these constants will be underlined. In addition, 
there will be constants for all those objects outside of $V$ which our forcing will add; those constants will be dotted.
\begin{align*} 
{\underline T} & \mbox{ \ \ \ intended to denote } T \\
{\underline x} \mbox{ for every } x \in H_\kappa  & \mbox{ \ \ \ intended to denote } x \mbox{ itself } \\ 
{\dot n} \mbox{ for every } n<\omega & \mbox{ \ \ \ as terms for elements of } M_i \mbox{ and } N_i, i < \omega_1 \\
{\dot M}_i \mbox{ for } i<\omega_1 & \mbox{ \ \ \ intended to denote } M_i \\
{\dot \pi}_{ij} \mbox{ for } i \leq j \leq \omega_1 & \mbox{ \ \ \ intended to denote } \pi_{ij} \\
{\dot {\vec M}} & \mbox{ \ \ \ intended to denote } (M_j , \pi_{jj'} \colon j \leq j' \leq 
\omega_1^{N_i}) \mbox{ for } i < \omega_1 \\
{\dot N}_i \mbox{ for } i < \omega_1 & \mbox{ \ \ \ intended to denote } N_i \\
{\dot \sigma}_{ij} \mbox{ for } i \leq j < \omega_1 & \mbox{ \ \ \ intended to denote } \sigma_{ij} \\
{\dot a} & \mbox{ \ \ \ intended to denote the distinguished $a$-predicate of } M_i, N_i, i<\omega_1 \\
{\dot I} 
& \mbox{ \ \ \ intended to denote the distinguished ideal of } N_i, i<\omega_1\\
{\dot X}_\delta \mbox{ for } \delta<\omega_1 & \mbox{ \ \ \ intended to denote } X_\delta.
\end{align*}
As $T \in H_{\kappa}$, the first line above is redundant. The constants ${\dot n}$, $n<\omega$, will produce the term models $N_i$ for $i<\omega_1$; it is of course not important to use (dotted) natural numbers as these constants, the elements of any other fixed countable set (in $H_\kappa$) would be equally good. 

The formulae of ${\cal L}$ will be exactly the expressions of the following form.\footnote{Again, the first order language associated with ${\mathbb P}_{\rm max}$-structures $(M;\in,I,a)$ is the language of set theory augmented by predicates for $I$ and $a$.}
\begin{gather*} 
\ulcorner {\dot N}_i \models \varphi({\underline \xi}_1, \ldots , {\underline \xi}_k , {\dot n}_1 , \ldots , {\dot n}_\ell, {\dot a}, {\dot I} , 
{\dot M}_{j_1}, \ldots, {\dot M}_{j_m} , {\dot \pi}_{q_1 r_1} , \ldots , {\dot \pi}_{q_s r_s} , {\dot {\vec M}}) \urcorner \\
\mbox{for } i<\omega_1,
\xi_1, \ldots , \xi_k < \omega_1, n_1, \ldots , n_\ell < \omega, j_1, \ldots, j_m < \omega_1,
q_1 \leq r_1 < \omega_1, \ldots, q_s \leq r_s < \omega_1 \\
\mbox{and for } \varphi \mbox{ being a formula of the first order language associated with ${\mathbb P}_{\rm max}$-structures,}
\end{gather*}
as well as:
\begin{align*}
\ulcorner {\dot \pi}_{i \omega_1}({\dot n})={\underline x} \urcorner  & \mbox{ \ \ \ for }
i < \omega_1 \mbox{ and } x \in H_{\omega_2} \\
\ulcorner {\dot \pi}_{\omega_1 \omega_1}({\underline x})={\underline x} \urcorner  & \mbox{ \ \ \ for }
x \in H_{\omega_2}\\ 
\ulcorner {\dot \sigma}_{ij}({\dot n})={\dot m} \urcorner  & \mbox{ \ \ \ for }
i \leq j < \omega_1, n, m<\omega \\
\ulcorner ({\underline {\vec u}},{\underline {\vec \alpha}}) \in {\underline T} \urcorner  & \mbox{ \ \ \ 
for } {\vec u} \in {}^{<\omega} \omega \mbox{ and } {\vec \alpha} 
\in {}^{<\omega} \omega_2 \\
\ulcorner {\underline \delta} \mapsto {\underline \mu} \urcorner  & \mbox{ \ \ \ for } \delta<\omega_1,
{\mu}<\kappa \\
\ulcorner {\underline x} \in {\dot X}_\delta \urcorner  & \mbox{ \ \ \ for } \delta<\omega_1, x \in H_\kappa
\end{align*}

Let us write ${\cal L}^\lambda$ for the collection of all ${\cal L}$-formulae except for the 
formulae which mention elements outside of $Q_\lambda$, i.e., except for the formulae
of the form $\ulcorner {\underline \delta} \mapsto {\underline \mu} \urcorner$
for $\delta < \omega_1$ and ${\lambda} \leq {\mu} < \kappa$ as well as $\ulcorner {\underline x} \in {\dot X}_\delta \urcorner$ for $\delta < \omega_1$ and $x
\in H_\kappa \setminus Q_\lambda$.
We may and shall assume that ${\cal L}$ is built in a canonical way so that 
${\cal L}^\lambda \subset Q_\lambda$ and in fact ${\cal L}^\lambda = {\cal L} \cap Q_\lambda$.

We say that a potential certificate 
\begin{eqnarray*}
{\mathfrak C} = \langle \langle M_i , \pi_{ij} , N_i , \sigma_{ij} \colon i \leq j \leq \omega_1 \rangle , 
\langle (k_n , \alpha_n ) \colon n<\omega \rangle ,
\langle \lambda_\delta , X_\delta \colon \delta \in K \rangle \rangle 
\end{eqnarray*}
as in
(\ref{certificate}) is {\em pre-certified} by a collection $\Sigma$ of ${\cal L}^\lambda$-formulae 
if and only if (C.1) through (C.8) are satisfied by ${\mathfrak C}$
and there are surjections\label{the_surjections} $e_i \colon \omega \rightarrow
N_i$ 
for $i<\omega_1$
such that the following hold true.
\begin{enumerate}
\item[($\Sigma$.1)] 
$\ulcorner {\dot N}_i \models \varphi({\underline \xi}_1, \ldots , {\underline \xi}_k , {\dot n}_1 , \ldots , {\dot n}_\ell, {\dot a}, {\dot I} , 
{\dot M}_{j_1}, \ldots, {\dot M}_{j_m} , {\dot \pi}_{q_1 r_1} , \ldots , {\dot \pi}_{q_s r_s} , {\dot {\vec M}})
\urcorner \in \Sigma$ iff 
\begin{enumerate} 
\item[(a)] $i<\omega_1$,
\item[(b)] $\xi_1$, $\ldots$, $\xi_k \leq \omega_1^{N_i}$, 
\item[(c)] $n_1$, $\ldots$ , $n_\ell < \omega$, 
\item[(d)] $j_1$, $\ldots$, $j_m \leq \omega_1^{N_i}$,
\item[(e)] $q_1 \leq r_1 \leq \omega_1^{N_i}$, $\ldots$, $q_s \leq r_s \leq \omega_1^{N_i}$, and 
\end{enumerate}
$$N_i \models \varphi(\xi_1, \ldots , \xi_k , e_i(n_1) , \ldots , e_i(n_\ell), A \cap \omega_1^{N_i}, I^{N_i} , M_{j_1}, \ldots
, M_{j_m}, \pi_{q_1 r_1}, \ldots , \pi_{q_s r_s}, {\vec M}){\rm , }$$ where $I^{N_i}$ is the distinguished ideal of $N_i$ and ${\vec M}=\langle M_j , \pi_{jj'} \colon
j \leq j' \leq \omega_1^{N_i})$,
%
%
\item[($\Sigma$.2)] $\ulcorner {\dot \pi}_{i \omega_1}({\dot n})={\underline x} \urcorner \in \Sigma$ iff $i < \omega_1$, $n<\omega$, and $\pi_{i \omega_1}(e_i(n))=x$,
\item[($\Sigma$.3)] $\ulcorner {\dot \pi}_{\omega_1 \omega_1}({\underline x})={\underline x} \urcorner \in \Sigma$ iff $x \in H_{\omega_2}$,
\item[($\Sigma$.4)] $\ulcorner {\dot \sigma}_{ij}({\dot n})={\dot m} \urcorner \in \Sigma$ iff $i \leq j < \omega_1$, $n$, $m<\omega$, and $\sigma_{ij}(e_i(n))=e_j(m)$,
\item[($\Sigma$.5)] letting $F$ with ${\rm dom}(F)=\omega$ be the monotone enumeration 
of the G\"odel numbers of all formulae of the form $$\ulcorner {\dot N}_0 \models \varphi({\dot n}_1, \ldots , {\dot n}_\ell, {\dot a}, {\dot I})
\urcorner$$ with $\ulcorner {\dot N}_0 \models \varphi({\dot n}_1, \ldots , {\dot n}_\ell, {\dot a}, {\dot I})
\urcorner \in \Sigma$, we have that $\ulcorner ({\underline {\vec u}},{\underline {\vec \alpha}}) \in {\underline T} \urcorner \in \Sigma$ iff
there is some $n<\omega$ such that $\langle {\vec u},{\vec \alpha} \rangle 
= \langle (F(m),\alpha_m ) \colon m< n \rangle$ and $F(m)=k_m$ for all $m<n$,\footnote{Here, $\langle (k_n,\alpha_n) \colon n<\omega \rangle$ is a component of ${\mathfrak C}$.} 
\item[($\Sigma$.6)] $\ulcorner {\underline \delta} \mapsto {\underline \mu} \urcorner \in \Sigma$ iff
$\delta \in K$ and ${\mu} = \lambda_\delta$, and
\item[($\Sigma$.7)] $\ulcorner {\underline x} \in {\dot X}_\delta \urcorner \in \Sigma$ iff
$\delta \in K$ and $x \in X_\delta$.
\end{enumerate}

We say that a potential certificate ${\mathfrak C}$ as in
(\ref{certificate}) is {\em certified} by a collection $\Sigma$ of formulae 
if and only if ${\mathfrak C}$ is pre-certified by $\Sigma$ and, in addition,
\begin{enumerate}
\item[($\Sigma$.8)]\label{(C.9)} if $\delta \in K$, then 
$[\Sigma]^{<\omega} \cap X_\delta \cap
E \not= \emptyset$ for every $E \subset {\mathbb P}_{\lambda_\delta}$ which is dense in ${\mathbb P}_{\lambda_\delta}$ and definable over the structure $$(Q_{\lambda_\delta}; \in , {\mathbb P}_{\lambda_\delta}  , A_{\lambda_\delta} )$$ from parameters in $X_\delta$.\footnote{Equivalently,
$[\Sigma]^{<\omega} \cap
E \not= \emptyset$ for every $E \subset {\mathbb P}_{\lambda_\delta} \cap 
X_\delta$ which is dense in ${\mathbb P}_{\lambda_\delta} \cap 
X_\delta$ and definable over the structure $$(X_\delta; \in , {\mathbb P}_{\lambda_\delta} \cap 
X_\delta , A_{\lambda_\delta} \cap X_\delta )$$ from parameters in $X_\delta$.}
\end{enumerate}

Item ($\Sigma$.8) is to play a crucial role in the proof that our forcing preserves stationary sets, see the proof of Lemma \ref{P_stat_pres}. If $p$ is a condition with $\ulcorner \underline{\delta} \mapsto \underline{\lambda} \urcorner \in p$, and if the predicate $A_\lambda$ guesses -- via ($\Diamond$) -- a name ${\dot C}$ for a club subset of $\omega_1$, then ($\Sigma$.8) will guarantee that $p \Vdash {\check \delta} \in {\dot C}$, see the proof of Claim \ref{claim2.10}.


\begin{definition}\label{semc} Let ${\mathfrak C}$ as in (\ref{certificate})
be a potential certificate relative to $A$, $T$, $H_{\omega_2}^V$, $\langle A_\mu\,:\,\mu\in C \cap \lambda\rangle$, and $\langle \mathbb P_\mu\,:\,\mu\in C\cap\lambda\rangle$. We call ${\mathfrak C}$ a {\em semantic
certificate relative to} $A$, $T$, $H_{\omega_2}^V$, $\langle A_\mu\,:\,\mu\in C\cap \lambda\rangle$, and $\langle \mathbb P_\mu\,:\,\mu\in C\cap\lambda\rangle$, or just a {\em semantic
certificate}, iff 
there is a 
collection $\Sigma$ of formulae such that ${\mathfrak C}$ is certified by $\Sigma$.
We call $\Sigma$ a {\em syntactic certificate}
iff there is a semantic certificate ${\mathfrak C}$ such that ${\mathfrak C}$ is certified by
$\Sigma$. 
\end{definition}

In the proofs of Lemma \ref{P_is_nonempty} and of Claim \ref{its-a-sem-c} we will run Definition \ref{semc} inside a generic iterate $M$ of $V$. The corresponding iteration map $\pi$ will move the parameters $A$, $T$, $H_{\omega_2}^V$, $\langle A_\mu\,:\,\mu\in C\cap \lambda\rangle$, and $\langle \mathbb P_\mu\,:\,\mu\in C\cap\lambda\rangle$, so that we will talk about semantic certificates relative to the parameters $\pi(A)$, $\pi(T)$, $\pi(H_{\omega_2}^V)$, $\pi(\langle A_\mu\,:\,\mu\in C \cap \lambda\rangle)$, and $\pi(\langle \mathbb P_\mu\,:\,\mu\in C\cap\lambda\rangle)$ with the understanding that this will be an object obtained by running Definition \ref{semc} over $M$ with these shifted parameters.  

Given a syntactic certificate $\Sigma$, there is a {\em unique} semantic
certificate ${\mathfrak C}$ such that ${\mathfrak C}$ is certified by $\Sigma$.
Even though it is obvious how to construct ${\mathfrak C}$ from $\Sigma$, in the proof of Lemma \ref{generics_give_certificates} below we \label{details} will provide details on how to derive a semantic certificate from a given $\Sigma$.

It is worth stressing that not every collection of ${\cal L}^\lambda$-formulae which is merely consistent is already a syntactic certificate. The requirement that the constant $x \in H_\kappa$ is to be interpreted by itself (cf.\ ($\Sigma$.2), ($\Sigma$.3), and
($\Sigma$.7)) may be restated as saying that for a consistent 
${\cal L}^\lambda$-theory to be a syntactic certificate it is to be true that certain types are omitted.

Let $\Sigma$ and $p$ be sets of formulae, where $p$ is finite.
We say that $p$
is {\em certified by} 
$\Sigma$ if and only if
there is some (unique) ${\mathfrak C}$ as in (\ref{certificate})
such that 
${\mathfrak C}$ is certified by 
$\Sigma$ and   
\begin{enumerate}
\item[($\Sigma$.9)] $p \in [\Sigma]^{<\omega}$.
\end{enumerate}
We may also say that $p$ is {\em certified by} ${\mathfrak C}$ as in (\ref{certificate}) iff there is some $\Sigma$ such that ${\mathfrak C}$ and $p$ are both certified by $\Sigma$---and we will then also refer to $\Sigma$ as a syntactical certificate for $p$ and to ${\mathfrak C}$ as the associated semantic certificate.

We are then ready to define the forcing ${\mathbb P}_\lambda$. 
We say that $p \in {\mathbb P}_\lambda$ if and only if 
\begin{eqnarray}\label{in_gen_ext}
V^{{\rm Col}(\omega,\lambda)} \models \mbox{``There is a set $\Sigma$ of ${\cal L}^\lambda$-formulae such that } p \mbox{ is certified by } \Sigma \mbox{.''}
\end{eqnarray}

Let $p$ be a finite set of formulae of ${\cal L}^\lambda$. By the homogeneity of ${\rm Col}(\omega,\lambda)$, if there is some $h$ which is 
${\rm Col}(\omega,\lambda)$-generic over $V$ and there is some 
$\Sigma
\in V[h]$
such that $p$ is certified by $\Sigma$, then for all $h$ which are ${\rm Col}(\omega,\lambda)$-generic over $V$ there is some $\Sigma \in V[h]$ such that $p$ is certified by $\Sigma$.
It is then easy to see that $\langle {\mathbb P}_\lambda \colon \lambda \in C \cup
\{ \kappa \} \rangle$ is definable over $V$ from $\langle A_\lambda \colon
\lambda < \kappa \rangle$
and $C$, and is hence an element of $V$.\footnote{To remind the reader, $C$ is the club from p.\ \pageref{defn-club}.} 

The following absoluteness fact will be relevant in the proofs of Lemma \ref{P_is_nonempty} and of Claim \ref{its-a-sem-c}.

\begin{lemma}\label{abs_of_ex_of_certificate}
Let $\lambda \in C \cup \{ \kappa \}$, and let $p$ be a finite set of formulae of ${\cal L}^\lambda$. If there is any outer model in which there is some $\Sigma$ which
certifies $p$, then there is some $\Sigma \in V^{{\rm Col}(\omega,\lambda)}$
which certifies $p$.
\end{lemma}

{\sc Proof.}
The statement that there is a $\Sigma$ which
certifies $p$ is $\Sigma_1$ in the parameters $A$, $T$, $H_{\omega_2}^V$, $\langle A_\mu \colon \mu \in C \cap \lambda \rangle$, and $\langle {\mathbb P}_\mu \colon \mu \in C \cap \lambda \rangle$, all of which are elements of ${\rm HC}^{V^{{\rm Col}(\omega,\lambda)}}$. Hence
by Shoenfield absoluteness (\ref{el_equivalent_shoenfield}), see \cite[Corollary 7.21]{book}, if there is any outer model in which there is some $\Sigma$ which
certifies $p$, then there is some $\Sigma \in V^{{\rm Col}(\omega,\lambda)}$
which certifies $p$.\footnote{In fact, if $P$ is a transitive model of {\sf KP} plus the axiom
{\em Beta} with $A$, $T$, $H_{\omega_2}^V$, $\langle A_\mu \colon \mu \in C \cap \lambda \rangle$, $\langle {\mathbb P}_\mu \colon \mu \in C \cap \lambda \rangle \in P$ and if $p \in {\mathbb P}_\lambda$,
then there is some $\Sigma \in P^{{\rm Col}(\omega,\lambda)}$ which certifies $p$.}\label{abs_of_cert}  \hfill $\square$ (Lemma \ref{abs_of_ex_of_certificate})

\subsection{Some properties of the forcing.}

It is easy to see that \begin{enumerate}
\item[(i)] ${\mathbb P}={\mathbb P}_\kappa \subset H_\kappa$,
\item[(ii)] if ${\bar \lambda} < \lambda$ are both in $C \cup \{ \kappa \}$, then ${\mathbb P}_{\bar \lambda} \subset {\mathbb P}_\lambda$, and
\item[(iii)] if $\lambda \in C \cup \{ \kappa \}$ is a limit point of 
$C \cup \{ \kappa \}$, then ${\mathbb P}_\lambda = \bigcup_{{\bar \lambda} \in C \cap \lambda} \, {\mathbb P}_{\bar \lambda}$,
\end{enumerate}
so that there is some club $D \subset C$ such that for all $\lambda \in D$,
$${\mathbb P}_\lambda = {\mathbb P} \cap Q_\lambda.$$ Hence ($\Diamond$) gives us the following.
\begin{enumerate}\label{diamond-thing}
\item[($\Diamond({\mathbb P})$)] For all $B \subset H_\kappa$ the set
$$\{ \lambda \in C \colon (Q_\lambda;\in , {\mathbb P}_\lambda , A_\lambda) \prec (H_\kappa;\in , {\mathbb P} , B) \}$$
is stationary.
\end{enumerate}

The first one of the following lemmas is entirely trivial.

\begin{lemma}\label{+++}
Let $\Sigma$ be a syntactic certificate, and let $p$, $q \in [\Sigma]^{<\omega}$. Then $p$ and $q$ are compatible conditions in ${\mathbb P}$. 
\end{lemma}

\begin{lemma}\label{P_is_nonempty}
$\emptyset \in {\mathbb P}_{\min(C)}$.
\end{lemma}

{\sc Proof.} This is a simple variant of the proofs of \cite[Theorem 2.8]{david-ralf} and of 
\cite[Theorem 4.2]{both}. What needs to be done is to construct a semantic/syntactic certificate (for $\emptyset$) in some outer model. 

Let
$h$ be ${\rm Col}(\omega,\omega_2)$-generic over $V$.
Let $q^* = (N^*;\in,I^*, a^*) \in ({\mathbb P}_{\rm max})^{V[h]}$
be as in the paragraph preceeding (\ref{q-star}), i.e., 
$q^* \in ({\mathbb P}_{\rm max})^{V[h]}$, $q^* < (H_{\omega_2}^V; \in,
{\sf NS}_{\omega_1}^V,A)$, $q^* \in D^*$, and such that (\ref{q-star}) is true, i.e., $q^* \in p[{T}] \subset p[{\tilde T}]$, cf.\ Lemma \ref{ex_of_T}.
Let $(M_i, \pi_{ij} \colon i \leq j \leq \omega_1^{N^*}) \in N^*$ be the unique generic iteration of the $({\mathbb P}_{\rm max})^{V[h]}$-condition $(H_{\omega_2}^V; \in,
{\sf NS}_{\omega_1}^V,A)$ which witnesses that $q^*$ is stronger than this condition. 

Let $(N_i,\sigma_{ij} \colon i \leq j \leq \kappa) \in V[h]$ be a generic iteration of
$N_0 = N^*$ such that $\kappa=\omega_1^{N_\kappa}$.\footnote{If we wished, we could even arrange that
writing $N_\kappa = (N_\kappa;\in,I',a')$, we have that $I' = ({\sf NS}_\kappa)^{V[h]}
\cap N_\kappa$, but this is not relevant here; cf.\ footnote \ref{no-requirement}.} Let 
\begin{eqnarray}\label{it_of_p}
(M_i , \pi_{ij} \colon i \leq j \leq \kappa)=
\sigma_{0 \kappa}((M_i, \pi_{ij} \colon i \leq j \leq \omega_1^{N_0}))
\end{eqnarray}
Since $M_{0} =  ((H_{\omega_2})^V; \in, ({\sf NS}_{\omega_1})^V, A)$ and $({\sf NS}_{\omega_1})^V$ is assumed to be saturated in $V$,
every maximal antichain in $V$ consisting of stationary subsets of $\omega_1$ is an element of $M_{0}$.
By \cite[Lemma 3.8]{hugh},
we may hence lift the generic ultrapower map  
$\pi_{01} \colon M_{0} \rightarrow M_{1}$ to act on all of $V$, and inductively we may  
lift the entire generic iteration (\ref{it_of_p})
to a generic iteration 
\begin{eqnarray}\label{it_of_pV}
(M_i^+ , \pi^+_{ij} \colon i \leq j \leq \kappa)
\end{eqnarray}
of $V$ in such a way that all $M_i^+$, $i \leq \kappa$, are transitive. 
Let us write $M=M^+_\kappa$ and $\pi = \pi^+_{0 \kappa}$.

Let $\langle k_n , \alpha_n \colon n<\omega \rangle$ be such that $x
= \langle k_n \colon n<\omega \rangle$ is
a real code for $N_0$ in the sense of Convention \ref{convention_real_code} and $\langle (k_n , \alpha_n) \colon n<\omega \rangle \in [T]$.
We then clearly have that $\langle (k_n , \pi(\alpha_n)) \colon n<\omega \rangle \in [\pi(T)]$.

\bigskip
\begin{tikzpicture}

\def\x{2.1};

\def\y{0.5};

\node (pT) at (0*\x,  0*\y){$p[T]$};

\node (ppiT) at (0*\x+1.5, 0*\y){$p[\pi(T)]$};

\node (N) at (0*\x, -2*\y) {$(N^*;\in, I^*, a^*)$};

\node (Nast) at (2*\x, -2*\y) {$N_\kappa$};

\node (M0) at (-2*\x, -4*\y) {$M_0$};

\node (MN) at (0*\x, -4*\y) {$M_{\omega_1^{N^*}}$};

\node (M3) at (2*\x, -4*\y) {$M_{\kappa}$};

\node (H) at (-2*\x, -6*\y) {$((H_{\omega_2})^V;\in, (\mathsf{NS}_{\omega_1})^V, A)$};

\node (V) at (-2*\x, -8*\y) {$V$};

\node (Mp) at (2*\x, -8*\y) {$M_\kappa^+$};

\node (M) at (2*\x+1, -8*\y) {$M$};

\path (pT)--(ppiT) node[midway]{$\subseteq$};

\path (N)--(pT) node[midway, rotate=90]{$\in$};

\draw[->] (N)--(Nast) node[midway, above]{$\sigma_{0\kappa}$};

\path (MN)--(N) node[midway, rotate=90]{$\in$};

\path (M3)--(Nast) node[midway, rotate=90]{$\in$};

\draw[->] (M0)--(MN) node[midway, above]{$\pi_{0 \omega_1^{N^*}}$};

\draw[->] (MN)--(M3) node[midway, above]{$\pi_{\omega_1^{N^*} \kappa}$};

\path (H)--(M0) node[midway, rotate=90]{$=$};

\path (M3)--(Mp) node[midway, rotate=270]{$\in$};

\path (H)--(V) node[midway, rotate=270]{$\in$};

\draw[->] (V)--(Mp) node[midway, above]{$\pi$};

\path (Mp)--(M) node[midway]{$=$};

\end{tikzpicture}

\bigskip
It is now easy to see that
\begin{eqnarray}\label{certificate_at_rho}
{\mathfrak C} = \langle \langle M_i , \pi_{ij} , N_i , \sigma_{ij} \colon i \leq j \leq \kappa \rangle , 
\langle (k_n , \pi(\alpha_n) ) \colon n<\omega \rangle ,
\langle \rangle \rangle
\end{eqnarray} 
certifies $\emptyset$ relative to the parameters $\pi(A)$, $\pi(T)$, $\pi(H_{\omega_2}^V)$, $\langle \rangle$, and $\langle \rangle$, 
with $\emptyset$ being construed as the empty set of $\pi({\cal L}^\kappa)$ formulae: as the third
component $\langle \rangle$ of ${\mathfrak C}$ in (\ref{certificate_at_rho}) is empty, {\em any} set of surjections $e_i \colon \omega \rightarrow N_i$, $i<\omega_1$, will induce a syntactic certificate for $\emptyset$ relative to $\pi(A)$, $\pi(T)$, $\pi(H_{\omega_2}^V)$, $\langle \rangle$, and $\langle \rangle$, whose associated semantic certificate is ${\mathfrak C}$. 
The statement that there is a syntactic certificate for $\emptyset$ is $\Sigma_1$ in the parameters $H_{\omega_2}^M$, $\pi(A)$, and $\pi(T)$,
which will all be in ${\rm HC}^{M^{{\rm Col}(\omega,\pi(\omega_2))}}$. Hence
by Shoenfield absoluteness (\ref{el_equivalent_shoenfield}), see \cite[Corollary 7.21]{book}, there is then some
${\mathfrak C} \in M^{{\rm Col}(\omega,\pi(\omega_2))}$ as in (\ref{certificate_at_rho}) which certifies $\emptyset$ relative to $\pi(A)$, $\pi(T)$, $\pi(H_{\omega_2}^V)$, $\langle \rangle$, and $\langle \rangle$, so that 
$\emptyset \in \pi(\mathbb P_{\min(C)})$, cf.\ Lemma \ref{abs_of_ex_of_certificate}.\footnote{Exactly in order to be able to do this we let the forcing also search for $\langle M_i , \pi_{ij} \colon i \leq j \leq \omega_1 \rangle$ rather than just $\langle N_i , \sigma_{ij} \colon i \leq j \leq \omega_1 \rangle$. The presence of
$\langle M_i , \pi_{ij} \colon i \leq j \leq \omega_1 \rangle$ allows us to lift $\pi_{0\kappa}$ to a map acting on all of $V$, so that we may then apply Shoenfield absoluteness and pull back the statement of interest---namely
$\emptyset \in \pi(\mathbb P_{\min(C)})$. Cf.\ the discussion on p.\ \pageref{discussion_why_M}.  \label{footnote}} By the elementarity of $\pi$,
then, 
$\emptyset \in {\mathbb P_{\min(C)}}$.
\hfill $\square$ (Lemma \ref{+++})

\begin{lemma}\label{generics_give_certificates}
Let $\lambda \in C \cup \{ \kappa \}$.
Let $g \subset {\mathbb P}_\lambda$ be a filter such that $g \cap E \not= \emptyset$ for all dense $E \subset {\mathbb P}_\lambda$ which are definable over $(Q_\lambda;\in,{\mathbb P}_\lambda)$ from elements of $Q_\lambda$. Then 
$\bigcup g$
is a syntactic certificate.
\end{lemma}

{\sc Proof.}
Let us first describe how to read off from $\bigcup g$ a candidate
$${\mathfrak C} = \langle \langle M_i , \pi_{ij} , N_i , \sigma_{ij} \colon i \leq j \leq \omega_1 \rangle , 
\langle (k_n , \alpha_n ) \colon n<\omega \rangle ,
\langle \lambda_\delta , X_\delta \colon \delta \in K \rangle \rangle$$
for
a semantic certificate for $\bigcup g$. 
A variant of what is to come shows 
how to derive ${\mathfrak C}$ from a given syntactic certificate $\Sigma$,
where ${\mathfrak C}$ is unique such that $\Sigma$ certifies ${\mathfrak C}$, cf.\ the remark on p.\ \pageref{details}.

Some of the formulas to follow simply describe the construction of a direct limit associated with the maps $\pi_{ij}$ and $\sigma_{ij}$. 

For 
$i$, $j <\omega_1$ and 
$\tau$, $\sigma \in \{ {\dot n} \colon n<\omega \} \cup \{ {\underline \xi} \colon \xi < \omega_1 \}$
define 
\begin{align*} 
\tau \sim_{i} \sigma & \mbox{ \ \ \ iff \ \ } \ulcorner {\dot N}_i \models \tau = \sigma \urcorner \in \bigcup g \\
(i,\tau) \sim_{\omega_1} (j,\sigma) & \mbox{ \ \ \ iff } i \leq j \wedge \exists
\rho \, \{ \ulcorner
{\dot \sigma}_{ij}(\tau)=\rho \urcorner , \ulcorner {\dot N}_j \models
\rho = \sigma \urcorner \} \subset \bigcup g\\
{} & \mbox{ \ \ \ \ or } j \leq i \wedge \exists
\rho \, \{ \ulcorner
{\dot \sigma}_{ji}(\sigma)=\rho \urcorner , \ulcorner {\dot N}_i \models
\rho = \tau \urcorner \} \subset \bigcup g\\
[\tau]_{i} & = \{ \sigma \colon \tau \sim_{i} \sigma \}  \\
[(i,\tau)] & = \{ (j,\sigma) \colon (i,\tau) \sim_{\omega_1} (j,\sigma) \}  \\
{M}_i & = \{ [\tau]_i \colon \tau \in \{ {\dot n} \colon n<\omega \} \cup \{ {\underline \xi} \colon \xi < \omega_1 \} \wedge \ulcorner {\dot N}_i \models \tau \in {\dot M}_i \urcorner \in \bigcup g \} \\
M_{\omega_1} & = (H_{\omega_2})^V \\
{N}_i & = \{ [\tau]_i \colon \tau \in \{ {\dot n} \colon n<\omega \} \cup \{ {\underline \xi} \colon \xi < \omega_1 \} \} \\
N_{\omega_1} & = \{ [i,\tau] \colon i<\omega_1 \wedge \ulcorner {\dot N}_i \models \tau = \tau \urcorner \in \bigcup g \} \\
[\tau]_{i} \, {\tilde \in}_i \, [\sigma]_{i} & \mbox{ \ \ \ iff \ \ } \ulcorner {\dot N}_i \models \tau \in \sigma \urcorner \in \bigcup g \\
[i,\tau] {\tilde \in}_{\omega_1} [j,\sigma] & \mbox{ \ \ \ iff }
i \leq j \wedge \exists
\rho \, \{ \ulcorner
{\dot \sigma}_{ij}(\tau)=\rho \urcorner , \ulcorner {\dot N}_j \models
\rho \in \sigma \urcorner \} \subset \bigcup g\\
{} & \mbox{ \ \ \ \ or } j \leq i \wedge \exists
\rho \, \{ \ulcorner
{\dot \sigma}_{ji}(\sigma)=\rho \urcorner , \ulcorner {\dot N}_i \models
\tau \in \rho \urcorner \} \subset \bigcup g\\
[\tau]_i \in I^{N_i} & \mbox{ \ \ \ iff \ \ } \ulcorner {\dot N}_i \models \tau \in {\dot I} \urcorner \in \bigcup g \\
[i,\tau] \in I^{N_{\omega_1}} & \mbox{ \ \ \ iff \ \ } [\tau]_i \in I^{N_i} \\
[\tau]_i \in a^{N_i} & \mbox{ \ \ \ iff \ \ } \ulcorner {\dot N}_i \models \tau \in {\dot a} \urcorner \in \bigcup g \\
[i,\tau] \in a^{N_{\omega_1}} & \mbox{ \ \ \ iff \ \ } [\tau]_i \in I^{N_i} \\ 
\pi_{ij}([\tau]_{i}) = [\sigma]_j & \mbox{ \ \ \ iff \ \ } \ulcorner {\dot N}_j \models 
{\dot \pi}_{ij}(\tau)=\sigma \urcorner \in \bigcup g \\
\pi_{i \omega_1}([\tau]_{i}) = x & \mbox{ \ \ \ iff \ \ } \ulcorner 
{\dot \pi}_{i \omega_1}(\tau)={\underline x} \urcorner \in \bigcup g \\
\pi_{\omega_1 \omega_1}(x) = x & \mbox{ \ \ \ iff \ \ } x \in (H_{\omega_2})^V \\
\sigma_{ij}([\tau]_{i}) = [\sigma]_j & \mbox{ \ \ \ iff \ \ } \ulcorner  
{\dot \sigma}_{ij}(\tau)=\sigma \urcorner \in \bigcup g \\
\sigma_{i \omega_1}([\tau]_i) = [i,\tau] & {} \\
%
(k,\alpha) = (k_n, \alpha_n) & \mbox{ \ \ \ iff } \exists {\vec u} \exists
{\vec \alpha} ( \ulcorner ({\underline {\vec u}} , {\underline {\vec \alpha}}) \in {\underline T} \urcorner
\in \bigcup g \wedge k = {\vec u}(n) \wedge  \alpha = {\vec \alpha}(n) )\\
\delta \in K & \mbox{ \ \ \ iff } \exists {\mu} \, \ulcorner {\underline \delta} \mapsto {\underline \mu} \urcorner \in \bigcup g\\
{\mu} = \lambda_\delta & \mbox{ \ \ \ iff } \delta \in K \wedge \ulcorner {\underline \delta} \mapsto {\underline \mu} \urcorner \in \bigcup g\\
x \in X_\delta & \mbox{ \ \ \ iff } \delta \in K \wedge \ulcorner {\underline x} \in {\dot X}_\delta \urcorner \in \bigcup g
\end{align*}
%

In order to see that this all works out we have to run a few density arguments. To show that a given subset of ${\mathbb P}$ is dense we frequently 
make use of Lemma \ref{+++}.
We will provide more details in some cases and fewer in others, and we are confident that the reader will be easily able to fill in the straightforward details herself in the latter cases.

Let us first observe that 
${\tilde \in}_0$ is wellfounded and that in fact (the transitive collapse of)
the structure $N_0 = (N_0;{\tilde \in}_0,a^{N_0},I^{N_0})$ is an iterable ${\mathbb P}_{\rm max}$ condition. This is true because of the following.

\begin{claim}\label{N_0}
(C.2) is true, i.e., $\langle (k_n, \alpha_n) \colon n<\omega \rangle \in [T]$ and $\langle k_n \colon
n<\omega \rangle$ codes the theory of $N_0$ in the sense of Convention \ref{convention_real_code}.
\end{claim}

{\sc Proof} of Claim \ref{N_0}. Let $m<\omega$. Writing $$q_0 = \{ \ulcorner (\underline{(k_n \colon n<m)},\underline{(\alpha_n \colon n<m)}) \in {\dot T} \urcorner \} {\rm , }$$ we have that $q_0 \in g.$ If 
$${\mathfrak C} = \langle \langle M'_i , \pi'_{ij} , N'_i , \sigma'_{ij} \colon i \leq j \leq \omega_1 \rangle , 
\langle (k'_n , \alpha'_n ) \colon n<\omega \rangle ,
\langle \lambda'_\delta , X'_\delta \colon \delta \in K' \rangle \rangle$$ certifies 
$q_0$, then $k'_n = k_n$ and $\alpha'_n = \alpha_n$ for $n<m$ by ($\Sigma.5$),
and then $$((k_n \colon n<m),(\alpha_n \colon n<m)) \in T$$ by (C.2). 

This shows $\langle (k_n, \alpha_n) \colon n<\omega \rangle \in [T]$.

By ($\Sigma.5$) and (C.2), for each $k<\omega$ the sets
\begin{eqnarray*}
D^0_k & = & \{ p \in {\mathbb P} \colon \exists m \, \exists ((k_n \colon n<m),(\alpha_n \colon n<m)) \exists r \, \\
{} & {} & ( \, \ulcorner (\underline{(k_n \colon n<m)},\underline{(\alpha_n \colon n<m)}) \in {\underline T} \urcorner \in p \wedge k_r = k \\
{} & {} &  \wedge \, k \mbox{ is the G\"odel number of } \ulcorner {\dot N}_0 \models \varphi({\dot n}_1 , \ldots , {\dot n}_\ell , {\dot a} , {\dot I}) \urcorner \, ) \rightarrow \\
{} & {} & \ulcorner {\dot N}_0 \models \varphi({\dot n}_1 , \ldots , {\dot n}_\ell , {\dot a} , {\dot I}) \urcorner \in p \}
\end{eqnarray*} and
\begin{eqnarray*}
D^1_k & = & \{ p \in {\mathbb P} \colon k \mbox{ is the G\"odel number of } \ulcorner {\dot N}_0 \models \varphi({\dot n}_1 , \ldots , {\dot n}_\ell , {\dot a} , {\dot I}) \urcorner \wedge \\
{} & {} & \ulcorner {\dot N}_0 \models \varphi({\dot n}_1 , \ldots , {\dot n}_\ell , {\dot a} , {\dot I}) \urcorner \in p \rightarrow \\
{} & {} & \exists m \, \exists ((k_n \colon n<m),(\alpha_n \colon n<m)) \exists r \, 
\, ( \, k_r = k \wedge \\
{} & {} & \ulcorner (\underline{(k_n \colon n<m)},\underline{(\alpha_n \colon n<m)}) \in {\underline T} \urcorner \in p \, )  \}
\end{eqnarray*}
are dense in ${\mathbb P}$. This implies that $\langle k_n \colon
n<\omega \rangle$ codes the theory of $N_0$ in the sense of Convention \ref{convention_real_code}.


\hfill $\square$ (Claim \ref{N_0})

\bigskip

Another set of easy density arguments will give that $(N_i ,\sigma_{ij} \colon i \leq j \leq \omega_1)$
is a generic iteration of $N_0$, where we identify $N_i$ with the structure
$(N_i;{\tilde \in}_i,a^{N_i},I^{N_i})$. 
To verify this, let us first show:
\begin{claim}\label{claim2.4}
For each $i < \omega_1$ and for each $\xi \leq \omega_1^{N_i}$, $[\underline{\xi}]_i$ represents 
$\xi$ in (the transitive collapse of the well-founded part of) the term model for $N_i$;
moreover, $a^{N_i} = A \cap \omega_1^{N_i}$. Hence $a^{N_{\omega_1}}=A$.
\end{claim}

{\sc Proof} of Claim \ref{claim2.4}. The set $$D^2 = \{ p \in {\mathbb P} \colon
\exists \xi \, \ulcorner {\dot N}_i \models \underline{\xi} \mbox{ is the least uncountable cardinal } \urcorner \in p \}$$ is easily seen to be dense in ${\mathbb P}$, so that there is some (unique!) $\xi_0<\omega_1$ such that writing $$p_0 = \{ \ulcorner {\dot N}_i \models \underline{\xi}_0 \mbox{ is the least uncountable cardinal } \urcorner \} {\rm , }$$
$p_0 \in g$.
Let us now prove by induction on $\xi \leq \xi_0$ 
that $[\underline{\xi}]_i$ must always represent $\xi$ in (the transitive collapse of the well-founded part of) $N_i$. Fix such $\xi$. 

For all $n<\omega$,
$$D_{n}^3 = \{ p \in {\mathbb P} \colon \ulcorner {\dot N}_i \models {\dot n} \in 
\underline{\xi} \urcorner \in p \rightarrow \exists \underline{\zeta} < \xi \, \ulcorner {\dot N}_i \models {\dot n} = 
\zeta \urcorner \in p \}$$
is dense below $p_0$. Also, for all $\zeta < \xi$,
$$D_{\zeta}^4 = \{ p \in {\mathbb P} \colon \ulcorner {\dot N}_i \models 
\underline{\zeta} \in \underline{\xi} \urcorner \in p \}$$ is dense below $p_0$. This shows that if $\tau \in \{ {\dot n} \colon n<\omega \} \cup \{ \underline{\xi} \colon \xi \in \xi_0+1 \}$ is a term, then $[\tau]_i {\tilde \in}_i [\underline{\xi}]_i$ iff $[\tau]_i = [\underline{\zeta}]_i$ for some $\zeta<\xi$. Using the inductive hypothesis, this then implies that $[\underline{\xi}]_i$ represents $\xi$ in (the transitive collapse of the well-founded part of) $N_i$. In particular, $\xi_0 = \omega_1^{N_i}$.

Now if $\xi \in A \cap \omega_1^{N_i}$, then $$D^5_\xi = \{ p \in {\mathbb P} \colon \ulcorner {\dot N}_i \models \underline{\xi} \in {\dot a} \urcorner \}$$ is dense below $p_0$, and if $\xi \in \omega_1^{N_i} \setminus A$, then $$D^6_\xi = \{ p \in {\mathbb P} \colon \ulcorner {\dot N}_i \models \underline{\xi} \notin {\dot a} \urcorner \}$$ is dense below $p_0$.  
Claim \ref{claim2.4} then follows. \hfill $\square$ (Claim \ref{claim2.4})

\bigskip
Similarly:

\begin{claim}\label{claim3.5} Let $i<\omega_1$. 
$N_{i+1}$ is generated from ${\rm ran}(\sigma_{ii+1}) \cup \{ \omega_1^{N_i} \}$
in the sense that for every $x \in N_{i+1}$ there is some function $f \in {}^{\omega_1^{N_i}} 
(N_i) \cap N_i$ such that $x = \sigma_{ii+1}(f)(\omega_1^{N_i})$.
\end{claim}

{\sc Proof} of Claim \ref{claim3.5}. Let $p_0$ be as in the proof of Claim \ref{claim2.4}, let $p \leq p_0$, and let $\Sigma \supset p$ be a syntactic certificate for $p$ with associated semantic certificate 
$${\mathfrak C} = \langle \langle M'_i , \pi'_{ij} , N'_i , \sigma'_{ij} \colon i \leq j \leq \omega_1 \rangle , 
\langle (k'_n , \alpha'_n ) \colon n<\omega \rangle ,
\langle \lambda'_\delta , X'_\delta \colon \delta \in K' \rangle \rangle . $$
Fix $i<\omega_1$, and let $e_i \colon \omega \rightarrow N'_i$ and $e_{i+1} \colon \omega \rightarrow N'_{i+1}$ be as on p.\ \pageref{the_surjections}. Notice that $\omega_1^{N'_i} = \omega_1^{N_i}$.

Let $n<\omega$. There must be some $f \in {}^{\omega_1^{N_i}}N_i' \cap N_i'$ with $e_{i+1}(n) = \sigma_{i i+1}'(f)(\omega_1^{N_i})$. Let $m$, $m' <\omega$ be such that $f=e_i(m)$ and
$\sigma_{i i+1}'(f)=e_{i+1}(m')$. Then 
\begin{flalign*}
& p \cup \{ \ulcorner {\dot N}_{i+1} \models {\dot n} = {\dot m}'(\underline{\omega_1^{N_i}}) \urcorner , \ulcorner {\dot \sigma}_{i i+1}({\dot m})={\dot m}' \urcorner \} \leq p.
\end{flalign*}
This argument shows that the set 
$$D^7_n = \{ p \in {\mathbb P} \colon
\exists m \, \exists m' \, \{ \ulcorner {\dot N}_{i+1} \models {\dot n} = {\dot m}'(\underline{\omega_1^{N_i}}) \urcorner , \ulcorner {\dot \sigma}_{i i+1}({\dot m})={\dot m}' \urcorner \} \subset p \}$$ is dense below $p_0$.
Claim \ref{claim3.5} then follows.  
\hfill $\square$ (Claim \ref{claim3.5})

\begin{claim}\label{claim3.6}
Let $i<\omega_1$. $\{ X \in {\cal P}(\omega_1^{N_i}) \cap N_i \colon \omega_1^{N_i} \in 
\sigma_{ii+1}(X) \}$ is generic over $N_i$ for the forcing given by the $I^{N_i}$-positive sets.
\end{claim}

{\sc Proof} of Claim \ref{claim3.6}. Let $p_0$, $p$, $\Sigma$, ${\mathfrak C}$, $e_i$, and $e_{i+1}$ be as in the previous proof. Let $n<\omega$ be such that $e_i(n)$ is a maximal antichain in $N_i'$ for the forcing given by the $I^{N_i'}$-positive sets. Let $m$, $m' <\omega$ be such that $e_i(m) \in e_i(n)$ and $\omega_1^{N_i} = \omega_1^{N_i'} \in \sigma'_{i i+1}(e_i(m)) = e_{i+1}(m')$. Then 
\begin{flalign*}
& p \cup \{ \ulcorner {\dot N}_i \models {\dot m} \in {\dot n} \urcorner ,
\ulcorner {\dot \sigma}_{i i+1}({\dot m})={\dot m}' \urcorner , \ulcorner {\dot N}_{i+1} \models \underline{\omega_1^{N_i}} \in\dot m' \urcorner \} \leq p.
\end{flalign*}
This argument shows that the set 
$$D^8_n = \{ p \in {\mathbb P} \colon
\exists  m \, \exists m' \, \{ \ulcorner {\dot N}_i \models {\dot m} \in {\dot n} \urcorner ,
\ulcorner {\dot \sigma}_{i i+1}({\dot m})={\dot m}' \urcorner , \ulcorner {\dot N}_{i+1} \models \underline{\omega_1^{N_i}} \in {\dot m}' \urcorner \} \subset p \}$$
is dense below $p_0$.
Claim \ref{claim3.6} then follows.  
\hfill $\square$ (Claim \ref{claim3.6})

\bigskip
Claims \ref{claim3.5} and \ref{claim3.6} readily imply that if $i<\omega_1$, then $N_{i+1}$ is a generic ultrapower of $N_i$.
By the next claim, direct limits are taken at limit stages:

\begin{claim}\label{claim3.7}
Let $i \leq \omega_1$ be a limit ordinal. For every $x \in N_i$ there is some $j<i$ and some $z \in N_j$ such that $x = \sigma_{ji}(z)$.
\end{claim}

{\sc Proof} of Claim \ref{claim3.7}. This is trivial for $i=\omega_1$. Now let $i<\omega_1$. Let $p$, $\Sigma$, and ${\mathfrak C}$ be as in the previous two proofs. Fix a limit ordinal $i$. For each $n<\omega$ there are $j<i$ and $n' < \omega$ with $\sigma_{ji}(e_j(n'))=e_i(n)$, where $e_j \colon \omega \rightarrow N_j'$ and $e_i \colon \omega \rightarrow N_i'$ are as on p.\ \pageref{the_surjections}. Then 
\begin{flalign*}
& p \cup \{ \ulcorner {\dot \sigma}_{ji}({\dot n}')={\dot n} \urcorner \} \leq p.
\end{flalign*} 
This argument shows that the set 
$$D^9_n = \{ p \in {\mathbb P} \colon \exists n' \, \ulcorner {\dot \sigma}_{ji}({\dot n}')={\dot n} \urcorner \in p \}$$ is dense in ${\mathbb P}$.
Claim \ref{claim3.7} then follows.  
\hfill $\square$ (Claim \ref{claim3.7})

\bigskip
$(N_i ,\sigma_{ij} \colon i \leq j \leq \omega_1)$
is then indeed a generic iteration of $N_0$. As $N_0$ is iterable, we may and shall identify $N_i$ with its transitive collapse, so that 
(C.4) holds true.

Another round of density arguments will show that
${\mathfrak C}$
satisfies (C.1), (C.3), (C.5), (C.6), and (C.7), where we identify $M_i$
with the structure $(M_i;\in,({\sf NS}_{\omega_1^{M_i}})^{M_i},A \cap \omega_1^{M_i})$.
Let us now verify (C.8) and ($\Sigma$.8), without writing down the relevant dense sets any more.

As for (C.8), its second part, $X_\delta \cap \omega_1 = \delta$ for $\delta \in K$, is easy. 
We will now use the Tarski-Vaught test to verify the first part of (C.8). 
Let $\varphi$ be any formula, and let $x_1$, $\ldots$, $x_k 
\in X_\delta$, $\delta \in K$.
Suppose that 
\begin{eqnarray}\label{in-the-model}
(Q_{\lambda_\delta};\in , {\mathbb P}_{{\lambda_\delta}},
A_{{\lambda_\delta}}) \models \exists v \, \varphi(v,x_1,\ldots , x_k).
\end{eqnarray}
Let $p \in g$ be such that $\{ \ulcorner {\underline x}_1 \in {\dot X}_\delta \urcorner , \ldots ,
\ulcorner \underline{x}_k \in {\dot X}_\delta \urcorner , \ulcorner \delta \mapsto \lambda_\delta \urcorner
\} \subset p$. Let $q \leq p$, and let
$\Sigma$ be a syntactical certificate for $q$ whose associated semantic certificate is
$${\mathfrak C}' =\langle \langle M'_i , \pi'_{ij} , N'_i , \sigma'_{ij} \colon i \leq j \leq \omega_1 \rangle , 
\langle (k'_n , \alpha'_n ) \colon n<\omega \rangle ,
\langle \lambda'_\delta , X'_\delta \colon \delta \in K' \rangle\rangle.$$
Then $\delta \in K'$ and $$\{ x_1 , \ldots , x_k \} \subset X'_{\delta} 
\prec (Q_{\lambda_\delta};\in , {\mathbb P}_{{\lambda_\delta}},
A_{{\lambda_\delta}}){\rm , }$$
so that by (\ref{in-the-model}) we may choose some $x \in X'_{\delta}$ with
$$(Q_{\lambda_\delta};\in , {\mathbb P}_{{\lambda_\delta}},
A_{{\lambda_\delta}}) \models \varphi(x,x_1,\ldots , x_k).$$ Let $r = q \cup \{
\ulcorner \underline{x} \in {\dot X}_{\delta} \urcorner \}$.

By density, there is then some $y \in X_\delta$ such that $$(Q_{\lambda_\delta};\in , {\mathbb P}_{{\lambda_\delta}},
A_{{\lambda_\delta}}) \models \varphi(y,x_1,\ldots , x_k).$$

The proof of ($\Sigma$.8) is similar.
Let again $\delta \in K$. Let $E \subset {\mathbb P}_{\lambda_\delta} \cap X_\delta^g$ be dense in ${\mathbb P}_{\lambda_\delta} \cap X_\delta$,
and $r \in E$ iff $r \in {\mathbb P}_{\lambda_\delta} \cap X_\delta$
and
\begin{eqnarray}\label{defn_of_E}
(Q_{\lambda_\delta};\in , {\mathbb P}_{{\lambda_\delta}},
A_{{\lambda_\delta}}) \models \varphi(r,x_1, \ldots , x_k).
\end{eqnarray}
Let $p \in g$ be such that $\{ \ulcorner {\underline x}_1 \in {\dot X}_\delta \urcorner , \ldots ,
\ulcorner {\underline x}_k \in {\dot X}_\delta \urcorner , \ulcorner \delta \mapsto \lambda_\delta \urcorner
\} \subset p$. Let $q \leq p$, and again let
$\Sigma$ be a syntactical certificate for $q$ whose associated semantic certificate is
$${\mathfrak C}' = \langle \langle M'_i , \pi'_{ij} , N'_i , \sigma'_{ij} \colon i \leq j \leq \omega_1 \rangle , 
\langle (k'_n , \alpha'_n ) \colon n<\omega \rangle ,
\langle \lambda'_\delta , X'_\delta \colon \delta \in K' \rangle \rangle.$$
Then $[\Sigma]^{<\omega} \cap X'_\delta$ has an element, say $r$, such that (\ref{defn_of_E}) holds true. Let $$s = q \cup r \cup \{ \ulcorner \underline{r} \in {\dot X}_\delta \urcorner \}$$

By density, then, $g \cap X_\delta \cap E \not= \emptyset$. \hfill $\square$ (Lemma \ref{generics_give_certificates})

\bigskip
Forcing with any ${\mathbb P}_\lambda$ makes $\omega_2^V$ $\omega$-cofinal, as the iteration map $\pi_{0 \omega_1}$ as being added by the generic filter maps the ordinals of the countable model $M_0$ cofinally into $\omega_2^V$. If $\lambda < \kappa$ (and $\lambda \in C$), then ${\mathbb P}_\lambda$ has size $\aleph_2$, so that by a result of S.\ Shelah, see \cite[Corollary 23.20]{jech}, ${\mathbb P}_\lambda$ will collapse $\omega_1$ to become countable. We are now going to prove that ${\mathbb P} = {\mathbb P}_\kappa$, on the other hand, does not collapse $\omega_1$ and in fact preserves stationary subsets of $\omega_1$.

\subsection{The forcing preserves stationary sets.}\label{preservation-subsection}

\begin{lemma}\label{P_stat_pres}
Let $g$ be ${\mathbb P}$-generic over $V$. Let 
$${\mathfrak C} = \langle \langle M_i , \pi_{ij} , N_i , \sigma_{ij} \colon i \leq j \leq \omega_1 \rangle , 
\langle (k_n , \alpha_n ) \colon n<\omega \rangle ,
\langle \lambda_\delta , X_\delta \colon \delta \in K \rangle\rangle$$
be the semantic certificate associated with the syntactic certificate $\bigcup g$.
Let 
$$N_{\omega_1} = (N_{\omega_1};\in,A,I^*).$$ Then every element of
$({\cal P}(\omega_1) \cap N_{\omega_1}) \setminus I^*$ is stationary in $V[g]$.
\end{lemma}

\begin{corollary}\label{stat_sets_are_preserved} 
${\mathbb P}$ preserves stationary subsets of $\omega_1$.
\end{corollary}

{\sc Proof} of Corollary \ref{stat_sets_are_preserved} from Lemma \ref{P_stat_pres}.
Let ${\mathfrak C}$ 
be as in the statement of Lemma \ref{P_stat_pres},
and let us write $M_i = (M_i;\in,I_i,a_i)$ and $N_i = (N_i;\in,I_i^*,a_i^*)$ for $i \leq \omega_1$.
In the light of Lemma \ref{generics_give_certificates}, by (C.3) we will have that $I_{\omega_1^{N_0}} = I_0^* \cap M_{\omega_1^{N_0}}$, so that also $I_{\omega_1} = I^*_{\omega_1} \cap M_{\omega_1}$. By   
(C.5), the universe of $M_{\omega_1}$ is $(H_{\omega_2})^V$ and 
$I_{\omega_1}=({\sf NS}_{\omega_1})^V$, while $I_{\omega_1}^*$ is denoted by $I^*$ in the statement of Lemma \ref{P_stat_pres}. We thus get that $({\sf NS}_{\omega_1})^V = I^* \cap V$, so that the conclusion of Lemma \ref{P_stat_pres} also gives that ${\mathbb P}$ preserves stationary subsets of $\omega_1$. \hfill $\square$(Corollary \ref{stat_sets_are_preserved})

\bigskip
{\sc Proof} of Lemma \ref{P_stat_pres}. Let ${\dot N}_{\omega_1} \in
V^{\mathbb P}$
be a canonical name for $N_{\omega_1}$, and let ${\dot I}^* \in
V^{\mathbb P}$
be a canonical name for $I^*$.
Let ${\bar p} \in g$, ${\dot C}$, ${\dot S} \in V^{\mathbb P}$, and
$i_0<\omega_1$ and $n_0<\omega$ be such that
\begin{enumerate}
\item[(i)] ${\bar p} \Vdash$ ``${\dot C} \subset \omega_1$ is club,'' 
\item[(ii)] ${\bar p} \Vdash$ ``${\dot S} \in ({\cal P}(\omega_1) \cap {\dot N}_{\omega_1}) \setminus {\dot I}^*$,'' and \label{ithroughiii}
\item[(iii)] ${\bar p} \Vdash$ ``${\dot S}$ is  
represented by $[i_0,{\dot n_0}]$ in the term model producing ${\dot N}_{\omega_1}$.'' 
\end{enumerate}
We may and shall also assume that
\begin{eqnarray}\label{right-thing}
\ulcorner {\dot N}_{i_0} \models {\dot n_0} \mbox{ is a subset of the first uncountable cardinal, yet } {\dot n_0} \notin {\dot I} \urcorner \in {\bar p}{\rm , }
\end{eqnarray} 
because the ${\cal L}$-formula in (\ref{right-thing}) must belong to every syntactic certificate for ${\bar p}$, as ${\bar p}$ satisfies (ii) and (iii).

Let $p \leq {\bar p}$ be arbitrary, $p \in {\mathbb P}$. We aim to produce some $q \leq p$ and some $\delta<\omega_1$ such that $q \Vdash {\check \delta} \in {\dot C} \cap {\dot S}$,
see Claim \ref{claim2.10} below. 

For $\xi<\omega_1$, let $$D_\xi = \{ q \leq p \colon \exists \eta \geq \xi \, (\eta < \omega_1 \wedge q \Vdash {\check \eta} \in {\dot C}) \}{\rm, }$$ so that $D_\xi$ is open dense below $p$. Let $$E = \{ (q,\eta) \in {\mathbb P} \times \omega_1 \colon
q \Vdash {\check \eta} \in {\dot C} \}.$$ Let us write $$\tau = ((D_\xi \colon \xi<\omega_1),E).$$
We may and shall identify $\tau$ with some subset of $H_\kappa$ which codes $\tau$. Here and in what follows, $\kappa$ is still equal to $\omega_3$.

By ($\Diamond({\mathbb P})$), we may pick some 
$\lambda \in C$ such that $p \in {\mathbb P}_\lambda$ and
\begin{eqnarray}\label{choice_of_lambda}
(Q_\lambda;\in,{\mathbb P}_\lambda,A_\lambda) \prec (H_\kappa;\in,{\mathbb P},\tau).
\end{eqnarray}

Let $h$ be ${\rm Col}(\omega,\omega_2)$-generic over $V$, and let 
$g' \in V[h]$ be a filter on ${\mathbb P}_\lambda$ such that 
$p \in g'$
and $g'$ meets every dense set which is definable over $(Q_\lambda;\in,{\mathbb P}_\lambda,A_\lambda)$ from parameters in $Q_\lambda$.
By Lemma \ref{generics_give_certificates}, $\bigcup g'$ is a syntactic certificate for $p$,
and we may let
$$\langle \langle M_i' , \pi_{ij}' , N_i' , \sigma_{ij}' \colon i \leq j \leq \omega_1 \rangle , 
\langle (k_n' , \alpha_n' ) \colon n<\omega \rangle ,
\langle \lambda'_\delta , X_\delta' \colon \delta \in K' \rangle \rangle$$
be the associated semantic certificate. In particular, $K' \subset \lambda$.

Let $S$ denote the subset of $\omega_1$ which is represented by $[i_0,{\dot n_0}]$ in the term model giving $N_{\omega_1}'$, so that 
if $N_{\omega_1}'=(N_{\omega_1}',\in,A,I')$, then by (\ref{right-thing}),
\begin{eqnarray}\label{S_positive}
S \in ({\cal P}(\omega_1) \cap N'_{\omega_1})
\setminus I'.
\end{eqnarray}
Notice that $\omega_1^{V[h]} = \omega_3^V = \kappa$. 
Inside $V[h]$, we may extend $\langle N'_i , \sigma'_{ij} \colon i \leq j \leq \omega_1 \rangle$ to a generic iteration $$\langle N'_i , \sigma'_{ij} \colon i \leq j \leq \kappa \rangle$$ such that 
\begin{eqnarray}
\omega_1 \in \sigma'_{\omega_1,\omega_1+1}(S).
\end{eqnarray}
This is possible as $\omega_1^{N'_{\omega_1}} = {\rm sup} \{ \omega_1^{N_j} \colon j<\omega_1 \} = \omega_1$ and by (\ref{S_positive}).
Let $$\langle M'_i , \pi'_{ij} \colon i \leq j \leq \kappa \rangle = \sigma_{0,\kappa}(\langle M'_i , \pi'_{ij} \colon i \leq j \leq \omega_1^{N'_0} \rangle){\rm , }$$ so that $\langle M'_i , \pi'_{ij} \colon i \leq j \leq \kappa \rangle$ is an extension of $\langle M'_i , \pi'_{ij} \colon i \leq j \leq \omega_1 \rangle$.

Since $M'_{\omega_1} =  ((H_{\omega_2})^V; \in, ({\sf NS}_{\omega_1})^V, A)$, cf.\ (\ref{Momega1}), and $({\sf NS}_{\omega_1})^V$ is assumed to be saturated in $V$,
every maximal antichain in $V$ consisting of stationary subsets of $\omega_1$ is an element of $M'_{\omega_1}$.
By \cite[Lemma 3.8]{hugh}, we may hence lift the generic ultrapower map  
$\pi'_{\omega_1 \omega_1+1} \colon M'_{\omega_1} \rightarrow M_{\omega_1+1}'$ to act on all of $V$, and inductively we may  
lift the entire generic iteration
$\langle M'_i , \pi'_{ij} \colon \omega_1 \leq i \leq j \leq \kappa \rangle$ to a generic iteration $$\langle M_i^+ , \pi_{ij}^+ \colon \omega_1 \leq i \leq j \leq \kappa \rangle$$ of $V$ with all $M_i^+$, $\omega_1 \leq i \leq \kappa$, being transitive. 
Let us write $M=M_\kappa^+$ and $\pi=\pi_{\omega_1,\kappa}^+$. 

\bigskip
\begin{tikzpicture}

\def\x{1.7};

\def\y{0.6};

\def\p{'};

\node (pT) at (0*\x,  0*\y){$p[T]$};

\node (ppiT) at (0*\x+1.5, 0*\y) {$p[\pi(T)]$};

\node (N0) at (0*\x, -2*\y) {$N_0\p$};

\node (N1) at (2*\x, -2*\y) {$N_{\omega_1}\p$};

\node (Nrho) at (4*\x, -2*\y) {$N_\kappa\p$};

\node (M0) at (-2*\x, -4*\y) {$M_0\p$};

\node (MN0) at (0*\x, -4*\y) {$M_{\omega_1^{N_0}}\p$};

\node (M1) at (2*\x, -4*\y) {$M_{\omega_1}\p$};

\node (Mrho) at (4*\x, -4*\y) {$M_{\kappa}\p$};

\node (H) at (2*\x, -6*\y) {$((H_{\omega_2})^V;\in, (\mathsf{NS}_{\omega_1})^V, A)$};

\node (V) at (2*\x, -8*\y) {$V$};

\node (Mrhop) at (4*\x, -8*\y) {$M_\kappa^+$};

\node (M) at (4*\x+1, -8*\y) {$M$};

\path (pT)--(ppiT) node[midway]{$\subseteq$};

\path (N0)--(pT) node[midway, rotate=90]{$\in$};

\draw[->] (N0)--(N1) node[midway, above]{$\sigma_{0\omega_1}\p$};

\draw[->] (N1)--(Nrho) node[midway, above]{$\sigma_{\omega_1\kappa}\p$};

\path (MN0)--(N0) node[midway, rotate=90]{$\in$};

\path (M1)--(N1) node[midway, rotate=90]{$\in$};

\path (Mrho)--(Nrho) node[midway, rotate=90]{$\in$};

\draw[->] (M0)--(MN0) node[midway, above]{$\pi_{0 \omega_1^{N_0\p}}\p$};

\draw[->] (MN0)--(M1) node[midway, above]{$\pi_{\omega_1^{N_0\p} \omega_1}\p$};

\draw[->] (M1)--(Mrho) node[midway, above]{$\pi_{\omega_1 \kappa}\p$};

\path (H)--(M1) node[midway, rotate=90]{$=$};

\path (Mrho)--(Mrhop) node[midway, rotate=270]{$\in$};

\path (H)--(V) node[midway, rotate=270]{$\in$};

\draw[->] (V)--(Mrhop) node[midway, above]{$\pi$};

\path (Mrhop)--(M) node[midway]{$=$};

\end{tikzpicture}

\bigskip
The key point is now that $\langle M_i' , \pi'_{ij} , N'_i , \sigma'_{ij} \colon i \leq j \leq \kappa \rangle$ may be used to extend $\pi \mbox{''} \bigcup g'$ to a
syntactic certificate 
\begin{eqnarray}\label{supset}
\Sigma \supset \pi \mbox{''} \bigcup g'
\end{eqnarray} 
for $\pi(p)$ 
in the following manner.
Let 
$K^* = K' \cup \{ \omega_1 \}$. For $\delta \in K'$, let $\lambda_\delta^* =
\pi(\lambda_\delta')$ and $X_\delta^* =
\pi \mbox{''} X_\delta'$. Also, write $\lambda_{\omega_1}^* = \pi(\lambda)$ and $X_{\omega_1}^* = \pi \mbox{''} Q_\lambda$. Notice that $\omega_1 \in
\pi(C)$, so that $K^* \subset \pi(C)$. 
Let 
$${\mathfrak C}^* =
\langle \langle M'_i , \pi'_{ij} , N'_i , \sigma'_{ij} \colon i \leq j \leq \kappa \rangle \mbox{, } 
\langle (k'_n , \pi(\alpha'_n) ) \colon n<\omega \rangle  \mbox{, } 
\langle \lambda_\delta^* , X_\delta^* \colon \delta \in K^* \rangle \rangle.$$


\begin{claim}\label{its-a-sem-c} 
${\mathfrak C}^*$ is a semantic certificate
for $\pi(p)$ relative to the parameters $\pi(A)$, $\pi(T)$, $\pi(H_{\omega_2}^V)$, $\pi(\langle A_\mu\,:\,\mu\in C \cap \lambda\rangle)$, and $\pi(\langle \mathbb P_\mu\,:\,\mu\in C\cap\lambda\rangle)$.
\end{claim}

{\sc Proof} of Claim \ref{its-a-sem-c}: First notice that $$\langle (k'_n , \pi(\alpha'_n) ) \colon n<\omega \rangle \in 
[\pi(T)]$$ Next, if $\delta \in K'$, then $$X_\delta^* = \pi \mbox{''} X_\delta' \prec
(\pi({Q}_{\lambda'_\delta});\in,\pi({\mathbb P}_{\lambda'_\delta}) , \pi(A_{\lambda'_\delta})),$$
and 
$\pi \mbox{''} g' \cap X_\delta^* =
\pi \mbox{''} (g' \cap X_\delta')$; as $\bigcup g'$ is a syntactic certificate for $p$, we thus have that
$\pi \mbox{''} g' \cap X_\delta^* \cap E \not= \emptyset$ for every $E \subset
\pi({\mathbb P}_{\lambda_\delta'})$ which is dense in $\pi({\mathbb P}_{\lambda_\delta'})$ and definable over the structure $(\pi({Q}_{\lambda'_\delta}); \in , \pi({\mathbb P}_{\lambda_\delta'}) , \pi(A_{\lambda_\delta'}))$ from parameters in $X_\delta^*$. Finally, $X_{\omega_1}^*= \pi \mbox{''} Q_\lambda$ and the choice of $g'$ imply that $\pi \mbox{''} g' \cap X_{\omega_1}^* \cap E \not= \emptyset$ for every $E \subset
\pi({\mathbb P}_{\lambda})$ which is dense in $\pi({\mathbb P}_{\lambda})$ and definable over the structure $$(\pi(Q_\lambda); \in , \pi({\mathbb P}_{\lambda}) , \pi(A_{\lambda}))$$ from parameters in $X_{\omega_1}^*$. 
This buys us that ${\mathfrak C}^*$ is indeed a semantic certificate
for $\pi(p)$ as an element of $\pi({\mathbb P})$, and that therefore there is some syntactic certificate $\Sigma$ as in
(\ref{supset}), relative to $\pi(A)$, $\pi(T)$, $\pi(H_{\omega_2}^V)$, $\pi(\langle A_\mu\,:\,\mu\in C \cap \lambda\rangle)$, and $\pi(\langle \mathbb P_\mu\,:\,\mu\in C\cap \lambda\rangle)$, such that ${\mathfrak C}^*$ is certified by $\Sigma$.
\hfill $\square$ (Claim \ref{its-a-sem-c})

\bigskip
Now let $[{\dot m_0}]_{\omega_1+1}$ represent $\sigma'_{\omega_1 \omega_1+1}(S)$ in the term model
for $N_{\omega_1+1}'$ provided by $\Sigma$, so that\footnote{Here, ${\dot \sigma}_{i_0 \omega_1+1}$ and ${\dot N}_{\omega_1+1}$ are terms of the language associated with $\pi({\mathbb P}_\lambda)$, and $\ulcorner {\dot \sigma}_{i_0 \omega_1+1}({\dot n_0}) = {\dot m_0} \urcorner$ and $
\ulcorner {\dot N}_{\omega_1+1} \models \underline{\omega_1} \in {\dot m_0} \urcorner$ are formulae of that language.}  
$$\{ \ulcorner {\dot \sigma}_{i_0 \omega_1+1}({\dot n_0}) = {\dot m_0} \urcorner ,
\ulcorner {\dot N}_{\omega_1+1} \models \underline{\omega_1} \in {\dot m_0} \urcorner \}
\subset \Sigma{\rm , }$$
in other words,
\begin{eqnarray}\label{half_of_new_p}
\pi(p) \cup \{ \ulcorner {\dot \sigma}_{i_0 \omega_1+1}({\dot n_0}) = {\dot m_0} \urcorner ,
\ulcorner {\dot N}_{\omega_1+1} \models \underline{\omega_1} \in {\dot m_0} \urcorner \}
\mbox{ is certified by } \Sigma .
\end{eqnarray}

Let us now define 
\begin{eqnarray}\label{defn_q}
q^* = \pi(p) \cup \{
\ulcorner {\dot \sigma}_{i_0 \omega_1+1}({\dot n_0}) = {\dot m_0} \urcorner ,
\ulcorner {\dot N}_{\omega_1+1} \models \underline{\omega_1} \in {\dot m_0} \urcorner ,
\ulcorner \underline{\omega_1} \mapsto \underline{\pi(\lambda)} \urcorner \}.
\end{eqnarray}

In the light of Lemma \ref{abs_of_ex_of_certificate}, we thus established the following.

\begin{claim}\label{claim_1}
$q^* \in \pi({\mathbb P})$, as being certified by $\Sigma$.
\end{claim}

The elementarity of
$\pi \colon V \rightarrow M$ then gives some $\delta<\omega_1$ 
such that
\begin{eqnarray}
q = p \cup \{ \ulcorner {\dot \sigma}_{i_0 \delta+1}({\dot n_0}) = {\dot m_0} \urcorner ,
\ulcorner {\dot N}_{\delta+1} \models \underline{\delta} \in {\dot m_0} \urcorner , \ulcorner \underline{\delta} \mapsto \underline{\lambda} \urcorner \} \in {\mathbb P}.
\end{eqnarray}

\begin{claim}\label{claim_2}\label{claim2.10}
$q \Vdash {\check \delta} \in {\dot C} \cap {\dot S}$.
\end{claim}

{\sc Proof} of Claim \ref{claim_2}. $q \Vdash {\check \delta} \in {\dot S}$ readily follows from $$\{ \ulcorner {\dot \sigma}_{i_0 \delta+1}({\dot n_0}) = {\dot m_0} \urcorner ,
\ulcorner {\dot N}_{\delta+1} \models \underline{\delta} \in {\dot m_0} \urcorner \} \subset q,$$ the fact that ${\bar p} \geq p$ forces that ${\dot S}$ is 
represented by $[i_0,{\dot n_0}]$ in the term model giving ${\dot N}_{\omega_1}$, and the fact that
by Claim \ref{claim2.4}, $[\underline{\delta}]_{\delta+1}$ represents $\delta$ in the model $N_{\delta+1}$
of any semantic certificate for $q$ as being given by a generic which contains $q$. 

Let us now show that $q \Vdash {\check \delta} \in {\dot C}$. 
We will in fact show that $q$ forces that ${\check \delta}$ is a limit point of ${\dot C}$.
Otherwise
there is some $r \leq q$ and some $\eta < \delta$ such that 
\begin{eqnarray}\label{21}
r \Vdash {\dot C} \cap {\check \delta} \subset {\check \eta}.
\end{eqnarray}
Suppose that $r$ is certified by $\Sigma$, so that there is some
\begin{eqnarray}\label{prime}
\langle \langle M'_i , \pi'_{ij} , N'_i , \sigma'_{ij} \colon i \leq j \leq \omega_1 \rangle , 
\langle (k'_n , \alpha'_n ) \colon n<\omega \rangle ,
\langle \lambda'_{\bar \delta} , X'_{\bar \delta} \colon {\bar \delta} \in K' \rangle \rangle
\end{eqnarray}
which is certified by $\Sigma$ and 
$r \in [\Sigma]^{<\omega}$. We must have that 
\begin{enumerate}
\item[(a)] $\delta \in K'$, 
\item[(b)] $X'_\delta \prec (Q_\lambda; \in , {\mathbb P}_\lambda , A_\lambda)$,
\item[(c)] $X'_\delta \cap \omega_1 = \delta$,
and
\item[(d)] 
$[\Sigma]^{<\omega} \cap X'_\delta \cap
E \not= \emptyset$ for every $E \subset {\mathbb P}_{\lambda}$ which is dense in ${\mathbb P}_{\lambda} \cap 
X'_\delta$ and definable over the structure 
$$(Q_\lambda;\in, {\mathbb P}_\lambda , A_\lambda)$$
from parameters in $X'_\delta$.
\end{enumerate}
Here, (a) is given by $\ulcorner \underline{\delta} \mapsto \underline{\lambda} \urcorner \in r$,
(b) and (c) are given by (C.8), while (d) is exactly what ($\Sigma$.8) on p.\ \pageref{(C.9)} buys us. 

We have that $A_\lambda = \tau \cap Q_\lambda$, and hence $A_\lambda$ may be identified 
with the ordered pair $((D_\xi \cap Q_\lambda \colon \xi < \omega_1), E \cap Q_\lambda )$.
As $\eta < \delta \subset X_\delta'$, $D_\eta$ is definable over
the structure $$(Q_\lambda;\in, {\mathbb P}_\lambda , A_\lambda)$$
from a parameter in $X'_\delta$. By
(\ref{choice_of_lambda}), $D_\eta \cap Q_\lambda$ is dense in ${\mathbb P}_\lambda$.
By (d) above, there is then some
$s \in [\Sigma]^{<\omega} \cap X'_\delta \cap D_\eta$. 

By (\ref{choice_of_lambda}) again, the unique smallest $\eta' \geq \eta$ with $s \Vdash {\check \eta}' \in {\dot C}$ must be in $X'_\delta$, hence $\eta' < \delta$ by (c) above. 
But now $s$ is compatible with $r$, as they are both finite subsets of the very same $\Sigma$ which certifies them (cf.\ 
Lemma \ref{+++}).
We have reached a contradiction with (\ref{21}).
\hfill $\square$ (Claim \ref{claim_2})

\bigskip
Now ${\dot C}$, ${\dot S}$, and ${\bar p} \in g$ were such that (i) through (iii) on p.\ \pageref{ithroughiii} hold true. We showed that the set of all $q \leq {\bar p}$ with $q \Vdash {\dot C} \cap {\dot S} \not= \emptyset$ is dense. As ${\dot C}$ was arbitrary, this buys us that ${\dot S}^g$ will be stationary in $V[g]$. But then as ${\dot S}$ was arbitrary, this means that every element of $({\mathcal P}(\omega_1) \cap N_{\omega_1}) \setminus I^*$ will be stationary in $V[g]$. \hfill $\square$ (Lemma \ref{P_stat_pres})

\section{Open questions.}

Woodin \cite{hugh} also introduced the axiom $(*)^+$ as a strengthening of $(*)$. $(*)^+$ says that there is some pointclass $\Gamma \subset {\mathcal P}({\mathbb R})$ and some filter $g \subset {\mathbb P}_{\max}$ such that
\begin{enumerate}
\item[(1)] $L(\Gamma,{\mathbb R}) \models {\sf AD}^+$,\footnote{${\sf AD}^+$ is a natural strengthening of ${\sf AD}$ which was introduced by H.\ Woodin, see e.g.\ \cite[Definition 9.6]{hugh}.}
\item[(2)] $g$ is ${\mathbb P}_{\rm max}$-generic over $L(\Gamma,{\mathbb R})$,
and
\item[(3)] ${\mathcal P}({\mathbb R}) \subset L(\Gamma,{\mathbb R})[g]$.
\end{enumerate}
See \cite[p.\ 908]{hugh}. While the main
result of the current paper gives a new twist to the question if ${\sf MM}$ is compatible with $(*)^+$, see \cite[p.\ 923, Question (15) a)]{hugh}, it also leaves this question wide open. See \cite{hugh-preprint}.

There is a strengthening of ${\sf MM}^{++}$, isolated by Viale \cite{viale}, which has strong completeness properties modulo forcing similar to those of $(*)$. This is the axiom ${\sf MM}^{+++}$. It says that a class $\mathbb T$ of towers of ideals with certain nice structural properties is dense in the category of stationary set preserving forcings; in other words, for every stationary set preserving forcing ${\mathbb P}$ there is a tower $\mathcal T$ in $\mathbb T$ such that ${\mathbb P}$ completely embeds into $\mathcal T$ in such a way that the quotient forcing preserves stationary sets in $V^{\mathbb P}$. ${\sf MM}^{+++}$ implies ${\sf MM}^{++}$, if $\kappa$ is an almost super-huge cardinal, then there is a partial order ${\mathbb P} \subset V_\kappa$ which forces ${\sf MM}^{+++}$, and if there is a proper class of almost super-huge cardinals, then ${\sf MM}^{+++}$ is complete for the theory of  the $\omega_1$-Chang model\footnote{The $\omega_1$-Chang model is the $\subseteq$-minimal transitive model of ${\sf ZF}$ containing all ordinals and closed under $\omega_1$-sequences. It can be construed as $\bigcup_{\alpha\in \Ord}L([\a]^{\aleph_1})$ and it includes $L(\mathcal P(\omega_1))$ as a definable submodel.} with respect to stationary set preserving partial orders forcing ${\sf MM}^{+++}$.

Schindler \cite[Definition 2.10]{both} introduces ${\sf MM}^{*,++}$ as a strengthening of ${\sf MM}^{++}$ by relaxing ``forceable by a stationary set preserving forcing'' to ``honestly consistent'' in an appropriate formulation of ${\sf MM}^{++}$, see
\cite{both}.

It remains open if either of ${\sf MM}^{+++}$ or ${\sf MM}^{*,++}$ is really stronger than ${\sf MM}^{++}$.
While Viale's ${\sf MM}^{+++}$ is known to be consistent modulo
a super-huge cardinal, it is open at this point if ${\sf MM}^{*,++}$ is consistent at all relative to large cardinals.

\end{document}